%% file: main.tex
\RequirePackage{ifpdf}
\newif\iffinalversion
\finalversiontrue
\ifpdf
  \iffinalversion
    \documentclass[final,pdftex,a4paper,notitlepage,british,intlimits]{article}
  \else
    \documentclass[pdftex,a4paper,notitlepage,british,intlimits]{article}  
  \fi
\else
  \iffinalversion
    \documentclass[final,hypertex,a4paper,notitlepage,british,intlimits]{article}
  \else
    \documentclass[a4paper,notitlepage,british,intlimits]{article}  
  \fi
\fi

\usepackage{amscd}

\iffinalversion

\usepackage[usenoshowkeys,usenoborder,usenocolorlink,usenosourcespecials,usehyperref, useenumitem]{load_packages}
\else
\usepackage[useshowkeys,useborder,usenocolorlink,usenosourcespecials,usehyperref, useenumitem]{load_packages}
\fi

\usepackage{mycommands}

\author{Antongiulio Fornasiero%
\thanks{Universit\`a~di~Pisa,
Dipartimento~di~Matematica~``L.~Tonelli'',
Largo~Bruno~Pontecorvo~5,
56127~Pisa,~Italy}}
\title{O-minimal spectrum}

\ifpdf 
  \hypersetup{%
    pdfsubject=%
     {Properties of the spectra of definable sets in an o-minimal structure.},
    pdfkeywords=%
     {Spectrum, o-minimality, rational type},
  }
\fi

\allsectionsfont{\sffamily}

\begin{document}
\maketitle
\setlength{\emergencystretch}{1em}
\iffinalversion
\else
 \setlength{\overfullrule}{5.0pt}
 \setlength{\hfuzz}{7pt}
\fi

\begin{abstract}
\sf
\input{abstract.txt}
\end{abstract}

\indent\textit{Key words:} Spectrum, o-minimality

\tableofcontents
\include{introduction}

\input{preliminary}

\input{special}

\input{functions}

\input{compact}

\input{definable}

\input{amalgam}



\bibliography{my_bibliography}
\bibliographystyle{abbrv}

\end{document}

%% file: abstract.txt
Let $X$ be a definable sub-set of some o-minimal structure.
We study the spectrum of $X$, in relation with the definability of types.

%% file: introduction.tex
\section{Introduction}
Let $M$ be an o-minimal structure, expanding a group, and $A \subseteq M^n$ be a definable \subseT.
The spectrum of~$A$, denoted by~$\spc A$, as defined by Pillay~\cite{PILLAY:1988}, is the set of complete types over~$A$, with the topology with the following basis of open sets:
\[
\set{\spc U: U \subseteq A \et U \text{ open and definable}}.
\]
Note that $\spc A$ has also the Stone topology,
which is finer than the spectral one.
When speaking of topological notions about~$\spc A$,
we refer to the spectral topology, unless explicitly said otherwise.

$\spc A$ is a $T_0$ space, but not a $T_1$ space:
namely, points are not closed.
In this situation, there is the so-called specialization order:
$x$~is a specialization of $y$ (written $x \leq y$)
iff $\cl x \subseteq \cl y$.

\cite{PILLAY:1988} proved that $\spc A$ is a spectral space.
Namely, $\spc A$~has a basis of quasi-compact open sets stable under finite intersections, and every irreducible closed set is the closure of a unique point.
Coste and Carral~\cite{COSTE-CARRAL:1983} study spectral spaces in general, and the normal ones in particular (since $A$ is definably normal, $\spc A$~is normal).

We continue the study of the properties of~$\spc A$.
The main theme is the relationship between the specialization order~$\leq$, the Rudin-Keisler order $\RKleq$, and the dichotomy between rational and irrational types (\cf\ Theorem~\ref{THM:SCLOSED}, Lemma~\ref{LEM:RATIONAL-0}, Corollary~\ref{COR:MAX-RAT}).

In \S\ref{SEC:PRELIMINARY}, we list some results on o-minimal structures and on topological spaces, which will be used in later sections.

In \S\ref{SEC:SPECIAL}, we collect some basic results on the spectrum of~$A$.
In particular, we prove that, for every $x \in \spc A$, $\cl{x}$ is totally ordered by the specialization order (Lemma~\ref{LEM:CHAIN}).

Any definable function $f : A \to B$ induces a function $\spf : \spc A \to \spc B$, which is continuous if $f$ is.
We say that $y \RKleq x$ iff $y = \spf (x)$ for some definable function~$f$.
In \S\ref{SEC:MAPS}, we study the properties of~$\spf$ and of~$\RKleq$.
In particular, we prove
that, if $y \leq x$, then $x \RKleq x$ (Theorem~\ref{THM:REALISATION}).

In the remainder of the article, we assume that $M$ expands a field.
A type $x$ is called \sclosed iff $x$ is closed in some $\spc A$, with $A$ definably compact.
In \S\ref{SEC:COMPACT}, we study the definable compactifications of definable sets, and give some results about \sclosed types (Definition~\ref{DEF:SCLOSED} and Theorem~\ref{THM:SCLOSED-MAP}).

The main results of the article are in \S\S\ref{SEC:RATIONAL} and~\ref{SEC:AMALGAM}.
In \S\ref{SEC:RATIONAL}, we investigate the relationship between rational and \sclosed types (Theorem~\ref{THM:SCLOSED}).
We further analyze the relationship between $\leq$, $\RKleq$ and rationality
(Theorem~\ref{THM:MAX-RAT}, Lemma~\ref{LEM:STEP-1}, Lemma~\ref{LEM:RATIONAL-0}, and Corollary~\ref{COR:MAX-RAT}).

In \S\ref{SEC:AMALGAM}, we study the amalgam of a rational and a \tirrational extension of~$M$ (Theorem~\ref{THM:AMALGAM}).

I thank A.~Berarducci for his help on conceiving and writing this article, and for the numerous discussions on these and similar topics.

%% file: preliminary.tex
\section{Preliminaries about o-minimal structures}
\label{SEC:PRELIMINARY}
Let $M$ be an o-minimal structure, expanding an ordered group.
Let $A$ be a definable sub-set of $M^k$, for some $k \in \Nat$.

In the following, if $M$ expands an ordered field, then $1$ will denote the neutral element of the multiplication.
Otherwise, $1$~will be some fixed element of $M$ such that $1 > 0$.
Definable will mean ``definable with parameters from~$M$'', unless explicitly stated otherwise.
\begin{notation}
Let $X$ be a topological space.
For every $Y \subseteq X$ 
we shall denote by $\clt X Y$ (or simply by $\cl Y$ if $X$ is clear from the context) the topological closure of $Y$ in $X$.
The frontier of $Y$ is
\[
\fr_X Y := \clt X Y \setminus Y.
\]
\end{notation}
\begin{remark}
If $A\subseteq B\subseteq M^k$ are non-empty and definable, then $\dim (\fr_B A) < \dim A$.
\end{remark}
\begin{definizione}\label{DEF:COMPACT}
Let $X \subseteq M^k$.
We say that $X$ is definably compact, or \intro{\dcompact} for short, iff $X$ is definable, closed and bounded.
$X$~is \intro{locally \dcompact} iff $X$ is definable, and for every $x \in X$ there exists a \dcompact neighborhood of~$x$ (in~$X$).
\end{definizione}
The reader can skip the remainder of this section,
and refer back to it when needed.
\begin{lemma}
Let $Z$ be a topological space.
Let $C, U \subseteq Z$ such that $U$ is open (in~$Z$).
Then,
\[
U \cap \clt Z C = U \cap \clt Z{U \cap C} = \clt U {U \cap C}.
\]
\end{lemma}
\begin{proof}
The fact that $U \cap \clt Z C \supseteq U \cap \clt Z{U \cap C}$ is obvious.
For the opposite inclusion, let $x \in U \cap \clt Z C$. Hence, $x \in U$.
If, for contradiction, $x \notin \clt Z{U \cap C}$, then
\[
x \in \clt Z{C \setminus U} \subseteq \clt Z{Z \setminus U} = Z \setminus U,
\]
because $U$ is open, absurd.
\end{proof}
\begin{corollary}\label{COR:LOCDIM}
If $C, U \subseteq A$ are definable, and $U$ is open, then
\[
\dim (U \cap C) = \dim(U \cap \cl C).
\]
\end{corollary}
\begin{proof}
$U \cap \cl C = \clt U {U \cap C}$, and $\dim\Pa{\clt U {U \cap C}} = \dim (U \cap C)$.
\end{proof}
\begin{lemma}\label{LEM:LOCDIM}
Let $C_1, C_2$ be definable disjoint \subsets of~$A$, and $m := \max\Pa{dim(C_1), \dim(C_2)}$.
Then,
\[
\dim\Pa{\cl{C_1} \cap \cl{C_2}} < m.
\]
\end{lemma}
\begin{proof}
In fact,
\begin{align*}
\cl{C_1} \cap \cl{C_2} &= \Pa{(\cl{C_1} \cap \cl{C_2}) \setminus C_1} 
\cup \Pa{(\cl{C_1} \cap \cl{C_2}) \setminus C_2} \subseteq \\
&\subseteq \Pa{\cl{C_1} \setminus C_1} \cup \Pa{\cl{C_2} \setminus C_2}.
\end{align*}
Since $\dim\Pa{\cl{C_i} \setminus C_i} < \dim C_i \leq m$ for $i =1, 2$, we are done.
\end{proof}
\begin{lemma}\label{LEM:SMALL-NEIGH}
Let $C_1$, $C_2 \subseteq A$ be closed and definable, and $C_0 := C_1 \cap C_2$.
Then, there exist $V_i \subseteq A$ open and definable, $i = 1, 2$, such that
\begin{enumerate}
\item $C_i \setminus V_i = C_0$, $i =1,2$;
\item $V_1 \cap V_2 = \emptyset$;
\item $\cl{V_1} \cap \cl{V_2} \subseteq C_0$.
\end{enumerate}
\end{lemma}
\begin{proof}
Define
\[
V_1 := \set{a \in A: d(a,C_1) < {\textstyle\frac{1}{3}} d(a,C_2)},
\]
where $d$ is the Euclidean distance (or equivalently the max distance if $M$ has no field structure), and similarly for~$V_2$.

Alternative proof:
let
\[\begin{aligned}
U &:= A \setminus C_0,\\
C_i' &:= C_i \setminus C_0, \quad i =1 ,2.
\end{aligned}
\]
$U$~is definable, and therefore definably normal.
Moreover, $C_1'$ and $C_2'$ are disjoint open \subsets of~$U$.
Therefore, there exist $V_1, V_2 \subseteq U$ disjoint and open in~$U$ (and hence open in~$X$) such that $C_i' \subseteq V_i$.
\end{proof}
\begin{definizione}
Let $X$ be a topological space.
$X$~is a $T_5$ space (also called \intro{completely normal}) iff every sub-space of $X$ is~$T_4$.%
\footnote{We do not assume that $T_4$ implies Hausdorff.}
\end{definizione}
Note that any metric space is~$T_5$.
Note moreover that $T_4$ does not imply $T_5$.
\begin{remark}
Let $X$ be a topological space.
The following are equivalent:
\begin{enumerate}
\item $X$ is $T_5$.
\item For every $U \subseteq X$, if $U$ is open, then $U$ is~$T_4$.
\item Let $C_1, C_2 \subseteq X$ be closed, and $C_0 := C_1 \cap C_2$.
Then, there exist $V_1, V_2 \subseteq X$ open, and satisfying conditions 1--3 of Lemma~\ref{LEM:SMALL-NEIGH}.
\item For every $D_1, D_2 \subseteq X$, if $D_1 \cap \cl{D_2} = \cl{D_1} \cap D_2 = \emptyset$, then there exist $V_1, V_2$ disjoint open \subsets of~$X$, such that $D_i \subseteq V_i$.
\end{enumerate}
\end{remark}
\begin{lemma}\label{LEM:NEIGH-F}
Let $C \subseteq U \subseteq D \subseteq M^k$ be definable, such that $D$ is \dcompact and $U$ is open in~$D$.
Let $f: D \to M$ be a definable continuous function, such that $C = f^{-1}(0)$.
Then, there exists $\varepsilon \in M$ such that $0 < \varepsilon$ and $f^{-1}\Pa{\inter{-\varepsilon} {\varepsilon}} \subseteq U$.
\end{lemma}
\begin{proof}
Note that $C$ must be closed.
Assume, for contradiction, that for every $t > 0$ there exists $a_t \in D \setminus U$ such that $\abs{f(t)} \leq t$.
By definable choice, we can find $a_t$ as above that is a definable function of~$t$.
Since $D$ is definably compact, there exists $\lim_{t \to 0^+} a_t =: a$.
However, $a \in C \setminus U$, a contradiction.
\end{proof}
Therefore, if $D$ is \dcompact and $f: D \to M$ is definable and continuous, the family $\set{f^{-1}\Pa{\inter{-t}{t}}: t > 0}$ is a fundamental system of open neighborhoods of~$f^{-1}(0)$.
Note that this is not true if $D$ is not \dcompact.

A variant of the following lemma is~\cite[Lemma~1.1]{WILKIE:2003}.
\begin{lemma}\label{LEM:RETRACT}
Let $C \subseteq A \subseteq M^n$ be definable, such that $C$ is a cell.
Then, there exists a definable open neighborhood $V$ of $C$ (in $A$) such that $C$ is a retract of~$V$.
Namely, there exists a definable continuous map $\rho: V \to C$, such that the restriction of $\rho$ to $C$ is the identity.
\end{lemma}
\begin{proof}
First, we will do the case when $A = M^n$.
The proof is by induction on~$n$.
If $n = 0$, the conclusion is trivial.
Since $C$ is a cell, there exists a cell $D \subseteq M^{n-1}$
such that one of the following cases happens:
\begin{enumerate}
\item $C = \Gamma(f) := \set{(x,y) \in M^n: x \in D \et f(x) = y}$, where $f: D \to M$ is a continuous definable function;

\item $C = \cell{f}{g}{D} \set{(x,y) \in M^n: x \in D \et f(x) < y <g(x)}$, where $f,g : D \to M \sqcup\mset{\pm \infty}$ are continuous definable functions such that $f < g$.
\end{enumerate}
By inductive hypothesis, there exist a definable open neighborhood $W$ of $D$ and a retraction $\sigma: W \to D$.
In the first case, let $V := W \times M$, and define
\[
\rho(w,t) := \Pa{\sigma(w), f(w)}.
\]
In the second case, define
\[\begin{aligned}
F:\ & W \times M \to D \times M\\
   & (w,t) \mapsto (\sigma(w),t).
\end{aligned}\]
Let $V := F^{-1}(C)$, and $\rho$ be the restriction of $F$ to~$V$.
Since $C$ is open in $D \times M$, $V$ is an open neighborhood of~$C$.

When $A \neq M^n$, let $V'$ and $\rho': V' \to C$ be the open neighborhood of $C$ in $M^n$ and the retraction whose existence we proved above.
Define $V := V' \cap A$, and $\rho := \rho' \rest V$.
\end{proof}
\begin{example}
Note that in general $C$ is not a retract of all~$M^n$.
For instance, let $n = 1$ and $C = \inter{0}{1}$.
\end{example}
\begin{lemma}\label{LEM:RET-FIELD}
When $M$ expands a real closed field,
in Lemma~\ref{LEM:RETRACT}, we can weaken the hypothesis to 
$C$ locally closed in~$A$, instead of $C$ cell.
\end{lemma}
\begin{proof}
It is~\cite[Proposition~8.3.3]{DRIES:1998}.
\end{proof}
\begin{question}
What happens if we drop the condition that $M$ expands a real closed field in Lemma~\ref{LEM:RET-FIELD}?
\end{question}
\begin{lemma}\label{LEM:DENSE-OPEN}
Let $X$ be a topological spaces, and $Y$, $U$ be sub-spaces of $X$ such that $Y$ is dense and $U$ is open in $X$.
Then, $\clt X{Y \cap U} = \clt X{U}$.
\end{lemma}
\begin{proof}
It is obvious that $\cl{Y \cap U} \subseteq \cl{U}$.
For the opposite inclusion, let $b \in \cl{U}$, and $V$ be an open neighborhood of~$b$.
Since $U$ is open, $W := V \cap U$ is also open, and since $b \in \cl{U}$, $W \neq \emptyset$.
Therefore, $(Y \cap U) \cap V = W \cap Y \neq \emptyset$, because $U$ is dense in~$X$, and thus $b \in \cl{Y \cap U}$.
\end{proof}
\begin{lemma}\label{LEM:DENSE-COMPACT}
Let $A$ be definable, closed and locally \dcompact, and $B$ be definable, such that $A \subseteq B$ and $A$ is dense in~$B$.
Then, $A$ is open in~$B$.
\end{lemma}
\begin{proof}
Fix $a \in A$.
Let $K$ be a \dcompact neighborhood of~$a$.
Let $U \subseteq K$ be a definable open neighborhood of~$a$ (in~$A$).
Therefore, there exists $V \subseteq B$ open and definable such that and $U := V \cap A$.
Therefore,
\[
A \supseteq K = \clt B{K} \supseteq \clt B{U} = \clt B{V \cap A}.
\]
By Lemma~\ref{LEM:DENSE-OPEN}, $\clt B{V \cap A} = \clt B V$.
Therefore, $V \subseteq \clt B V \subseteq A$, and therefore $A$ is a neighborhood of~$a$.
\end{proof}
\begin{lemma}\label{LEM:CELL-PROJ}
Let $C \subseteq M^{k + h}$ be a cell, and $\pi: M^{k + h} \to M^k$ be the projection on the first $k$ coordinates.
Then, exactly one of the following 2 things happens:
either $\pi\rest C$ is injective, or $\dim \pi(C) < \dim C$.
\end{lemma}
\begin{lemma}\label{LEM:INJECTIVE}
Let $f: M^h \to M^k$ be definable.
Then, there exists a decomposition of $M^h$ into definable sets $\set{C_i: i \leq n}$,
such that for every $i\leq n$, either \mbox{$f\rest{C_i}$} is injective,
or \mbox{$\dim f(C_i) < h$}.%
\end{lemma}
\begin{proof}%
\footnote{Thanks to prof.~Berarducci for the proof.}
Let
\[
S := \set{\Pa{f(a),a} \in M^{k+h}: a \in M^h}.
\]
Decompose $M^{k+h}$ into cells, in a way compatible with~$S$.
Let \mbox{$D \subseteq S$} be a cell, \mbox{$C := f^{-1}(D)$}, and $\pi: M^{k+h} \to M^k$ the projection on the first $k$ coordinates.
Note that \mbox{$f(C) = \pi(D)$}, and hence $\dim D \leq h$.
Moreover, $f\rest C$ is injective iff $\pi \rest D$ is injective.
By Lemma~\ref{LEM:CELL-PROJ}, either $\dim\pi(D) < h$, or $\pi\rest D$ is injective.
\end{proof}
\begin{lemma}\label{LEM:FRONTIER}
Let $X$ be a topological space, $Z \subseteq X$ be connected, and $U \subseteq X$ be open.
If $Z \cap U \neq \emptyset$ and $Z \setminus U \neq \emptyset$, then $Z \cap \fr U \neq \emptyset$.
\end{lemma}
\begin{proof}
Assume that $Z \cap \fr U = \emptyset$.
Then,
\[
Z = (Z \cap U) \sqcup (Z \setminus U) =
\Pa{Z \cap \cl{U}} \sqcup \Pa{Z \setminus U},
\]
contradicting the fact that $Z$ is connected.
\end{proof}

%% file: special.tex
\section{Specialization}
\label{SEC:SPECIAL}
Let $M$ be an o-minimal structure, expanding an ordered group.
Let $A$ be a definable sub-set of $M^k$, for some $k \in \Nat$, and $X := \spc A$ be the spectrum of~$A$ (namely, the set of complete types over~$A$).

In this section we study the basic properties of~$X$.
Most of these results are well-known, at least in the case when $M$ expands a field~\cite{COSTE-CARRAL:1983}, \cite{PILLAY:1987}, \cite{PILLAY:1988},
and~\cite{JONES:PHD}.
\begin{definizione}
Let $C \subseteq X$.
$C$~is \intro{definable} iff $C$ is of the form~$\spc D$,
where $D \subseteq A$ is definable.
$C \subseteq X$ is \intro{\tdef} iff $C$ is closed in the Stone topology, or, equivalently, $C$~is an intersection of definable sets.
\end{definizione}
\begin{definizione}
The \intro{spectral topology} in $X$ is the topology generated by the sets of the form $\spc U$, where $U \subseteq A$ is open and definable;
cf.~\cite{PILLAY:1988}.
\end{definizione}
\begin{remark}
Since the Stone topology is Hausdorff, any finite set is~\tdef.
Any definable set is~\tdef.
Since the Stone topology is stronger than the spectral one, any closed set is \tdef.
Since the Stone topology is compact, the spectral one is quasi-compact.
If $C\subseteq X$ is \tdef, and $f: A \to B$ is definable, then $\spf(C)$~is also \tdef.
\end{remark}
\begin{remark}
$\spc{\ }$ is an injective morphism of Boolean algebrae between the definable subsets of $A$ and the subsets of~$X$.
Moreover, if $B \subseteq A$ is definable, then
\[
\spc{\cl B} = \cl{\spc B}.
\]
\end{remark}
\begin{example}
However, $\spc{\ }$ does \emph{not} preserve infinite unions or infinite intersections.
For instance, let $A := [0,1]$, $U_t := \inter{-t}{t} \subseteq A$.
Then, $\bigcap_{0 < t \in M} U_t = \mset{0}$, but
\[
\bigcap_{0 < t \in M} \spc{U_t} = \set{0^-, 0, 0^+},
\]
where $0^+$ is the cut $\lOpenInter{-\infty}{0},\inter{0}{+ \infty}$, and similarly for~$0^-$.
\end{example}
\begin{lemma}
Let $C \subseteq X$ be open.
Then, $C$ is quasi-compact iff $C$ is definable.
\end{lemma}
\begin{proof}
The ``if'' direction is trivial.
On the other hand, since $C$ is open,  $C = \bigcup_{i\in I} U_i$, where each $U_i$ is open and definable.
If $C$ is quasi-compact, then there exists $I_0 \subseteq I$ finite, such that $C = \bigcup_{i \in I_0} U_i$, and thus $C$ is definable.
\end{proof}
\begin{definizione}
Let $x \in X$.
The \intro{dimension} of $x$ is
\[
\dim x := \min(\dim C),
\]
where $C$ varies among the definable \subsets of $X$ containing~$x$.
Given $C \subseteq X$ definable, the \intro{local dimension} of $C$ at $x$ is
\[
\locdim x C := \min(\dim U \cap C),
\]
where $U$ varies among the definable \emph{open}
\subsets of $X$ containing~$x$.
Define $\locdim x C = -\infty$ iff there exists an definable open $U$ containing $x$ such $U \cap C = \emptyset$.
Note that $\locdim x C \geq 0$ iff $x \in \cl C$, where $\cl C$ is the closure of $C$ in the spectral topology.

We shall write $\locdimo x$ instead of $\locdim x A$.

We say that $\LocDim(C)$ is constantly equal to $n$ iff $\locdim x C = n$ for every $x \in \cl C$.
Given $D \subseteq A$ definable, we define $\locdim x D := \locdim x {\spc D}$.
\end{definizione}
\begin{lemma}
$\dim x$ is the dimension of $c/M$, for any $N \ni c \models x$.
\end{lemma}
\begin{remark}
For every $x \in X$,
\[
\dim x \leq \locdimo x \leq \dim A.
\]
For every $C \subseteq X$ definable, Corollary~\ref{COR:LOCDIM} implies that
\[
\locdim x C = \locdim x {\cl C}.
\]
If $D \subseteq A$ is a cell of dimension~$n$, then $\LocDim D$ is constantly equal to~$n$.
If $\LocDim C$ is constantly equal to~$n$, then $\LocDim\Pa{\cl C}$ is also constantly equal to~$n$.
\end{remark}
\begin{remark}
$\locdimo x$ is the minimum of the dimensions of the definable open \subsets containing~$x$.
\end{remark}
Intuitively, the local dimension of $x$ is the dimension of the ambient space $A$ in a neighborhood of~$x$.
On the other hand, the dimension of $x$ tells us ``how big'' a definable \subseT containing $x$ must be.
\begin{example}
$X$~is not $T_1$ in general (namely, not every point is closed).
For instance, let $A := M$, $x = 0^+$.
Then, $\cl x = \set{0,0^+}$.
\end{example}
\begin{definizione}
For every $x, y \in X$, we say that $x$ is a \intro{specialization} of~$y$ (and  $y$ is a \intro{generalization} of~$x$) iff \mbox{$x \in \cl y$}, where $\cl y$ is the  closure of~$\mset{y}$ in the spectral topology, and we  write $x \leq y$ for this.
We shall say that $x$ is a closed point iff \mbox{$\cl x = \mset x$}, namely iff $x$ is minimal \wrt the order~$\leq$.
\end{definizione}
\begin{lemma}\label{LEM:ORDER}
\begin{enumerate}
\item $\leq$ is a partial order on $X$.
\item If $x \leq y$, then $\dim x \leq \dim y \leq \locdimo x$.
\item If $x \leq y$ and $D \subseteq X$ is definable, then $\locdim x D \geq \locdim y D$.
\item If $x < y$ (namely, $x \leq y$ and $x \neq y$), then $\dim x < \dim y$.
\end{enumerate}
\end{lemma}
\begin{proof}
1) The transitivity of $\leq$ is obvious.
To prove that $\leq$ is a partial order, it is enough to show that if $x \leq y$ and $y \leq x$, then $x = y$.
The hypothesis is equivalent to $\cl x = \cl y$.
Assume for contradiction that $x \neq y$.
Let $B, C \subseteq X$ be definable subsets such that
\[\begin{aligned}
x \in & B \setminus  C \text{ and}\\
y \in & C \setminus  B.
\end{aligned}\]
\Wlog, we can assume that $\dim B = \dim x \geq \dim y = \dim C$, and $B$ and $C$ disjoint.
Therefore, by Lemma~\ref{LEM:LOCDIM},
\[
\dim \Pa{\cl B \cap \cl C} < \dim B.
\]
Thus, $x \notin {\cl B} \cap {\cl C}$, and so $x \notin {\cl C}$.
However, $\cl y \subseteq \cl C$, absurd.

2)
Let $U \subseteq X$ be definable and open, such that $x \in U$ and $\dim U = \locdimo x$.
Moreover, $y \in U$, and therefore $\dim y \leq \dim U = \locdimo x$.

Assume, for contradiction, that $x \leq y$, but $m := \dim x > \dim y =: n$.
Let $C \subseteq X$ be definable such that $y \in C$ and $\dim C = n$.
Moreover, $\dim \cl C = \dim C = n$, therefore we can assume that $C$ is closed.
Besides, $x \notin {C}$, because $\dim x > \dim C$.
However, ${C}$~is closed and $y \in {C}$, and therefore $\cl y \subseteq {C}$, absurd.

3) Let $U\subseteq X$ be a open and definable, such that $\dim (U \cap D) = \locdim x D$, and $x \in U$. Then, $y \in U$, and we are done.

4) Assume, for contradiction, that $x < y$, but $\dim x = \dim y =: m$.
Let $E \subseteq A$ be definable, such that $y \in \spc E$, and $\dim E = m$.
Let $A' := \cl E$.
Then, $\dim A' = m$, and $x, y \in \spc{A'}$.
Therefore, \wloG, we can assume that $\dim A = m$.
Since $<$ is irreflexive, $y \notin \cl x$, namely there exist $U \subseteq X$ open and definable such that $y \in U$, but $x \notin U$.
However, $x \in {\cl U} \setminus {U}$, because $x \in \cl y$.
Therefore, $\dim x < m$, absurd.
\end{proof}
\begin{lemma}\label{LEM:CHAIN}
For every $x \in X$, $\cl x$ is totally ordered by~$\leq$.
Moreover, $\card{\cl x} \leq 1 + \dim x \leq 1 + \locdimo x \leq 1 + \dim A$.
\end{lemma}
\begin{proof}
Assume, for contradiction, that there exist $y_1, y_2 \in \cl x$, such that $y_1 \nleq y_2$ and $y_2 \nleq y_1$.
Let $D_i \subseteq A$ be definable such that
\[ \dim D_i = \dim y_i =: m_i, \quad \text{and}
\quad
y_i \in \spc{D_i} \setminus \spc{D_{2-i}}, \quad i =1,2.
\]
\Wlog, we can assume that $m_2 \leq m_1$.
\begin{claim}
$y_2 \leq y_1$.
\end{claim}
Note that the claim contradicts $y_2 \nleq y_1$, which is absurd.

Assume not.
Then, there exists $U \subseteq A$ open and definable, such that $y_2 \in \spc{U}$, but $y_1 \notin \spc U$.
By substituting $D_1$ with $D_1 \cap U$, and $D_2$ with $D_2 \setminus U$, we can assume that $D_2 \subseteq U$, and $D_1 \subseteq A \setminus U$, and in particular that $D_1 \cap D_2 = \emptyset$.

Let $C_i := \cl{D_i}$, $i = 1, 2$, and $C_0 := C_1 \cap C_2$.
Note that $C_0 \subseteq \cl U \setminus U$.
Hence, $y_2 \notin \spc{C_0}$.
Besides, since $D_1 \cap D_2 = \emptyset$, by Lemma~\ref{LEM:LOCDIM} we have $\dim C_0 < m_1$, and therefore $y_1 \notin \spc{C_0}$.
Let $V_1$ and $V_2$ be as in Lemma~\ref{LEM:SMALL-NEIGH}.
Since $V_1 \cap V_2 = \emptyset$, $x$~cannot be both in $\spc{V_1}$ and in~$\spc{V_2}$.
Assume that $x \notin \spc{V_i}$.
Since $y_i \notin C_0$, we have $y_i \in \spc{V_i}$.
However, $y_i \in \cl x \subseteq X \setminus \spc{V_i}$, absurd.

Therefore, we have proved that $\leq$ is a total order on $\cl x$.
Using Lemma~\ref{LEM:ORDER}, we conclude that  $\card{\cl x} \leq 1 + \dim x$.
\end{proof}
\begin{example}
Let $x \in A$.
The set $\set{y \in X: y \geq x}$ in general is \emph{not} totally ordered by~$\leq$.
For instance, let $A = M$, $x = 0$, $y = 0^+$ and $y' = 0^-$.
Then, $x < y, y'$, but neither $y \leq y'$ nor $y' \leq y$.
\end{example}
\begin{example}
Both cases are possible: $\card{\cl x} = \dim x$ and $\card{\cl x} < \dim x$.
For instance, let $A = M$. 
Let $x = 0^+$ and $x' = + \infty$.
Then, $\dim x = \dim x' = 1$.
However, $\cl x = \set{0, 0^+}$, while $\cl{x'} = \mset{x'}$.
\end{example}
\begin{remark}
If $x, y \in X$ are such that neither $x \leq y$ nor $y \leq x$, then there exist disjoint open definable \subsets $U, V \subseteq X$ such that $x \in U$ and $y \in V$.
\end{remark}
\begin{definizione}
Let $\iota : A \to X$ be the natural embedding, sending $c \in A$ to the type $x(a) := \text{``}a = c\text{''}$.
\end{definizione}
\begin{remark}
$\iota$~is a continuous map.
\end{remark}
\begin{remark}
Let $x \in X$. Then, $\dim x = 0$ iff $x = \iota(a)$ for some $a \in A$.
In that case, $x$~is a closed point.
\end{remark}
\begin{example}
Not all closed point have dimension~$0$.
For instance, if $A=M$, then $+\infty$ is a closed point of dimension~$1$.
\end{example}
\begin{remark}
Let $x \in X$.
Define the following partial types:
\[\begin{aligned}
\Phi(a) & :=
\set{a \in U: U \subseteq A \text{ definable and open},\ x \in \spc U},\\
\Psi(a) & :=
\set{a \in C: C \subseteq A \text{ definable and closed},\ x \in \spc C}.
\end{aligned}
\]
Then,
\[\begin{aligned}
\spc \Phi &= \set{y \in X: x \leq y},\\
\spc \Psi &= \set{y \in X: x \geq y} = \cl x.
\end{aligned}\]
Equivalently,
\[\begin{aligned}
\set{y \in X: x \leq y} & = \bigcap
\set{U: U \subseteq X \text{ definable and open},\ x \in U}\\
\cl x &= \bigcap
\set{C:  C \subseteq X \text{ definable and closed},\ x \in C}.
\end{aligned}\]
\end{remark}
\begin{remark}
More in general, if $D \subseteq X$, then
\[
\cl D = \bigcap \set{C: C\subseteq X \text{ definable and closed},\
D \subseteq C}.
\]
Moreover, $D$ is closed iff
\[
D = \bigcap \set{C: C\subseteq X \text{ definable and closed},\
D \subseteq C}.
\]
\end{remark}
\begin{corollary}
Let $x \in X$ such that $\dim x = \locdimo x$.
Then,
\[
\mset{x} = \bigcap
\set{U: U \subseteq X \text{ definable and open},\ x \in U}
\]
\end{corollary}
\begin{lemma}\label{LEM:MAX-EXTENSION}
Let $x \in X$.
Let $m := \dim x$, and $n := \locdimo x$.
Then, there exists $y \in X$ such that $x \leq y$ and $\dim y = n$.
\end{lemma}
\begin{proof}
If $m = n$, take $y := x$.
If $m < n$, let
\[
\Phi(a) := \set{a \in U \setminus C:\ C, U \subseteq A \text{ definable,}\
x \in \spc U,\ \dim C < n,\ U \text{ open}}.
\]
By definition of local dimension, it is easy to see that $\Phi$ is consistent set of formulae.
Let $y \in \spc \Phi$.
By definition of~$\Phi$, $\dim y = n$ and $x \leq y$.
\end{proof}
\begin{remark}
Let $\dim A = n$, and $k \leq n$.
The set $X_k := \set{x \in X: \dim x \geq k}$ is of the form $\spc \Phi_k$ (and hence $X_k$ is \tdef),
where $\Phi_k(a)$ is the partial type
\[
\set{a \notin C: C \subset A \text{ definable},\ \dim C < k}.
\]
Moreover, both $\iota(A)$ and $X_k$ are dense in~$X$, and, if $k > 0$, they are disjoint.
Finally, $\iota(A)$ is open in the Stone topology, because $X \setminus \iota(A) = X_1$.
\end{remark}
\begin{proof}
$\iota(A)$ is dense, because it is dense in the Stone topology, which is stronger than the spectral one.
$X_k$ is dense by Lemma~\ref{LEM:MAX-EXTENSION}.
\end{proof}
\begin{example}
Let $x, y \in X$ with $x < y$ and $m := \dim x$ and $n := \dim y$.
It is not true in general that if $m < l < n$,
then there exists $z \in X$ such that $x < z < y$ and $\dim z = l$.
For instance, let $A' := M^2$, $x = (0,0)$, $z = (0^+,0)$ ($z$ \emph{is} a complete type) and $y \in A'$ such that $\dim y = 2$ and $z < y$ ($y$~exists by Lemma~\ref{LEM:MAX-EXTENSION}).
Note that $x < z < y$.
Let $A$ be $A' := \set{(a,0): a > 0}$.
Since $z \in X' \setminus X$, and in $\cl y$ there is at most one $z'$ such that $\dim z' = 1$, there is no $z' \in X$ such that $z' < y$ and $\dim z' = 1$.
\end{example}
Let $A$ be definably compact.
Let $x, y \in X$ with $x < y$ and $m := \dim x$ and $n := \dim y$, and let
$l \in \Nat$ be such that $m < l < n$.
Later (examples~\ref{EX:GAP-EXP} and~\ref{EX:GAP-E}) we will show that
there might \emph{not} exist
$z \in X$ such that $x < z < y$ and $\dim z = l$.
\begin{lemma}\label{LEM:CLOSURE}
Let $C \subseteq X$ be \tdef.
Then,
\[
\cl C = \bigcup \set{\cl x: x \in C}.
\]
\end{lemma}
\begin{proof}
It is obvious that $\cl C \supseteq \bigcup\set{\cl x: x \in C}$.
For the other inclusion, we have to prove that for every $y \in \cl C$ there exists $x \in C$ such that $y \in \cl x$.
Since $C $ is \tdef, we can write $C = \bigcap_{i\in I} C_i$, where each $C_i$ is definable.
Let $\Phi$ be the following partial type
\[
\Phi(a) = \set{a \in U \cap C_i: i\in I,\
U \subseteq A \text{ open and definable},\ y \in \spc U}.
\]
Since $y \in \cl C$, $\Phi$ is consistent.
Any $x \in \spc \Phi$ satisfies $y \leq x$ and $x \in C$.
\end{proof}
\begin{corollary}
Let $C \subseteq X$.
Then, $C$ is closed iff $C$ is \tdef and 
\be\label{EQ:CLOSURE}
\forall x \in C\quad \cl x \subseteq C.
\ee
\end{corollary}
\begin{proof}
The ``only if'' direction is trivial.
For the ``if'' direction, Lemma~\ref{LEM:CLOSURE} implies that $\cl C = C$.
\end{proof}
\begin{example}
There are some \subsets $C \subset X$ that are not closed, but do satisfy \eqref{EQ:CLOSURE}.
By the Corollary, any such $C$ cannot be \tdef.
For instance, let $C := \iota(A)$.
$C$~is not closed (unless $A$ is finite), because $\cl C = X$, but $C$ does satisfy~\eqref{EQ:CLOSURE}.
\end{example}
\begin{lemma}
$A$~is definably connected iff $X$ is connected.
\end{lemma}
\begin{proof}
Easy.
\end{proof}
\begin{lemma}\label{LEM:THETA}
Let $M \elemeq N$, and
\[
\theta: \spc{A(N)} \to \spc{A(M)}
\]
be the restriction map.
Then, $\theta$ is continuous both in the Stone and the spectral topologies.
\end{lemma}
\begin{example}
$\theta$ is neither closed, nor open.
For instance, let $A = M$, $x = 0^+$, $c\in N$ such that $c \models 0^+$, $y := \iota (c) \in \spc{N}$.
Then, $y$~is a closed point, but $\theta(y) = x$ is not.
Let $U := \spc{\inter{-c}{c}} \subseteq \spc N$.
Then, $U$ is open, but $\theta(U) = \set{0^-,0,0^+}$ is not open.
\end{example}
\begin{question}
Is $\theta$ closed or open in the Stone topology?
\end{question}
\begin{lemma}\label{LEM:TAU}
Let $A, B$ be definable, and $\pi: A \times B \to A$ be the projection onto the first coordinate.
Let $N \extendeq M$, $c \in A(N)$, and $x := \tp c A \in X = \spc A$.
Let $\tau : \spc{B\Pa{M(c)}} \to \spc{A \times B}$ defined by
\[
\tau \Pa{y(b)} := x(a) \et \set{\phi(a,b): \phi(c,b) \in y(b)}.
\]
Then, $\tau$ is well-defined, and it is a (surjective) homeomorphism.
\end{lemma}
The case $B = [0,1]$ of the above Lemma is in~\cite{DELFS:1985}.
\begin{lemma}\label{LEM:CONNECTED}
Let $Z \subseteq X$ be \tdef and connected.
Let $\theta$ be as in the Lemma above.
Then, $\theta^{-1}(Z)$ is \tdef and connected.
\footnote{Thanks to Berarducci for the proof.}
\end{lemma}
\begin{proof}
Let $Z = \bigcap_{i \in I}{\spc{C_i}}$.
Then,
\[
W := \theta^{-1}(Z) = \bigcap_{i \in I}{\spc{C_i(N)}}.
\]
Assume, for contradiction, that $W$ is disconnected.
Namely, there exist $T_1$ and $T_2$ open such that
\begin{align*}
W \subseteq T_1 \cup T_2,\\
T_1 \cap T_2 \cap W = \emptyset,\\
W \cap T_i \neq \emptyset,&\quad i=1,2.
\end{align*}
\begin{claim}
We can also assume that the $T_i$ are definable (in~$N$).
\end{claim}
In fact,
\[
T_i = \bigcup_{j \in J_i} U_{i,j},
\]
where the $U_{i,j}$ are open and definable (in~$N$).
By compactness of the Stone topology, there exist $J'_i \subseteq J_i$ finite, $i=1,2$, such that
\begin{align*}
W \subseteq T_1' \cup T'_2 ,\\
W \cap T'_i \neq \emptyset,&\quad i=1,2,
\end{align*}
where $T'_i := \bigcup_{j \in J'_i} U_{i,j}$.
Moreover, $T'_1 \cap T'_2 \cap W \subseteq T_1 \cap T_2 \cap W = \emptyset$.

In particular, $\bigcap_{i \in I}{\spc{C_i(N)}} \subseteq T_1 \cup T_2$.
Using again the compactness of the Stone topology, we conclude that there exists $I_0 \subseteq I$ finite, such that
\begin{align*}
\spc{C(N)} \subseteq T_1 \cup T_2,\\
\spc{C(N)} \cap T_1 \cap T_2 = \emptyset,
\end{align*}
where $C := \bigcap_{i \in I_0}{C_i}$.
We have that $Z \subseteq \spc C$.
Since $Z$ is connected, $Z \subseteq \spc D$ for a (unique) definably connected component $D$ of $C$.
Hence, $W \subseteq \spc{D(N)} \subseteq T_1 \cup T_2$.
Moreover, $\spc{D(N)} \cap T_i \supseteq W \cap T_i \neq \emptyset$.
However, $D(N)$ is definably connected, hence $\spc{D(N)}$ is connected, a contradiction.
\end{proof}
\begin{conjecture}
Let $M \elemeq N$, $Z \subseteq X$ be \tdef, and $W := \theta^{-1}(Z)$,
where $\theta$ is as in Lemma~\ref{LEM:THETA}.
Then, $\theta$ induces an isomorphism between the \v{C}eck cohomology of~$Z$ and the one of~$W$.
\end{conjecture}
The conjecture is true if $M$ expands a field.
In fact, in that case we know that it holds if $Z$ is definable,
and therefore
\[
\Hc(Z) = \dirlim_{C \in \Deft {A} (Z)} \Hc(\spc C)
= \dirlim_{C \in \Deft {A}(Z)} \Hc\Pa{\theta^{-1}(\spc C)}
= \dirlim_{D \in \Deft {A(N)} (\theta^{-1}(Z))} \Hc(\spc D)
= \Hc\Pa{\theta^{-1}(\spc Z)},
\]
where
\[
\Deft{A}(Z) :=
\set{C \subseteq A: Z \subseteq C \et C \text{ is definable}}.
\]
The fact that $\Hc(Z) = \dirlim_{C \in \Deft {A} (Z)} \Hc(\spc C)$
will be proved elsewhere.
Note that the same proof works when $N$ is an o-minimal expansion of~$M$,
instead of an elementary extension.
\subsection{Beyond o-minimality}
Let $M$ be a first order topological structure (in the sense of Pillay~\cite{PILLAY:1987}).
We shall say that $M$ is definably $T_5$ iff every definable \subseT of $M^k$ is definably~$T_4$, for every $k \in \Nat$.
\begin{lemma}
The following are equivalent:
\begin{enumerate}
\item $M$~is definably $T_5$.
\item For every $U \subseteq M^k$, if $U$ is definable and open, then $U$ is definably~$T_4$.
\item  Lemma~\ref{LEM:SMALL-NEIGH} is true for any definable
$A \subseteq M^k$.
\item For every $D_1, D_2 \subseteq X$ definable, if $D_1 \cap \cl{D_2} = \cl{D_1} \cap D_2 = \emptyset$, then there exist $V_1, V_2$ disjoint definable open \subsets of~$X$, such that $D_i \subseteq V_i$.
\item For every $A \subseteq M^k$ definable, $\spc A$ is $T_4$.
\end{enumerate}
\end{lemma}
\begin{conjecture}
$M$ is definably $T_5$ iff, for every $A \subseteq M^k$ definable,
$\spc A$ is~$T_5$.
\end{conjecture}
Most of the results in this section apply to the following situation, with the same proofs (in particular, Lemmata~\ref{LEM:ORDER} and~\ref{LEM:CHAIN} hold).
$M$~is a first order  such that for every $k \in \Nat$ there is a function
\[
\dim : \Deft{M^k} \to [-1,k],
\]
where
\[
\Deft{M^k} := \set{A \subseteq M^k: A \text{ definable}},
\]
satisfying the following conditions:
\begin{align}
M &\text{ is definably } T_5;\\
\dim A &= -1 \quad\text{iff}\ \  A = \emptyset;\\
\dim(A \cup B) &= \max(\dim A, \dim B);\\
\dim(\fr A) &< \dim A \quad\text{if } A \neq \emptyset;\\
\dim\mset{a} &= 0\ \forall a \in M^k.
\end{align}
Note that $M$ must be~$T_2$.
In fact, $\dim\Pa{\cl a \setminus \mset a} = -1$ for every $a \in M^k$.
Hence, $M^k$ is~$T_1$, and, since $M$ is also definably~$T_4$, we have that $M$ is $T_2$ too.
\pagebreak[2]
\begin{remark}
If $A \subseteq M^k$ is definable, then $\spc A$ is a $T_4$ spectral space.
\end{remark}
\begin{proof}
The proof of \cite[Lemma~1.1]{PILLAY:1988} works also in this context.
More precisely, 
the definable open sets form a basis of quasi-compact open \subsets of~$\spc A$, stable under finite intersections.
Therefore, we need only to show that every irreducible closed set is the closure of a unique point.
Let $C \subseteq \spc A$ be closed and irreducible.
Note that $C$ \tdef, and hence compact in the Stone topology.
Let
\[
\Dfam := \set{D \subseteq \spc A: D \text{ is definable and closed } \et C \setminus D \neq \emptyset}.
\]
\begin{claim}
There exists $x \in C \setminus \bigcup \Dfam$.
\end{claim}
If, for contradiction, $C \subseteq \bigcup \Dfam$, then $\Dfam$ is a covering of $C$ by definable sets, and therefore, by compactness, there exists $D_1, \dotsc D_n \in \Dfam$ such that $C \subseteq \bigcup_{i = 1}^n D_i$.
Since $C$ is irreducible, $C \subseteq D_i$ for some $i \leq n$, absurd.

It is now easy to see that $C = \cl x$.
Uniqueness is a consequence of Lemma~\ref{LEM:ORDER}.

Since $A$ is definably~$T_5$, $\spc A$ is~$T_4$ (but we do not know whether it is~$T_5$).
\end{proof}
\begin{remark}
If $A \subseteq M^k$ is definable, then $A$ is a boolean combination of open definable sets.
\end{remark}
\begin{proof}
Induction on $\dim A$.
Since $\dim(\fr A) < \dim(A)$, we have that $\fr A$ is a boolean combination of open definable sets.
The conclusion follows from $A = \cl{A} \setminus \fr{A}$.
\end{proof}
Therefore, all the results in~\cite{COSTE-CARRAL:1983} about normal spectral spaces are true in this context.
Moreover, by Lemma~\ref{LEM:CHAIN}, the Krull dimension of $\spc A$ \cite[\P1.4]{COSTE-CARRAL:1983} is less or equal to~$\dim A$.
\begin{remark}
Let $A \subseteq M^k$ be definable.
If $A$ ha empty interior, then $\dim A < k$.
\end{remark}
\begin{proof}
Let $B := M^k \setminus A$.
Since $A$ has empty interior, $A = \fr B$.
\end{proof}

%% file: functions.tex
\section{Functions}
\label{SEC:MAPS}
Let $B \subseteq M^{h}$, for some $h \in \Nat$, $Y := \spc B$, and $f: A \to B$ be a definable function.
\begin{definizione}
Define
\[
\spc{f} : X \to Y
\]
as
\[
\spc f(x)(b) := \set{\phi(b): \phi\Pa{f(a)} \in x(a)}.%
\footnote{Here we use $a$ as a mute variable ranging in~$A$, and similarly for~$b$.}
\]
Namely,
\[
\spf(x) =
\bigcap\set{\spc{f(U)}: U \subseteq A \text{ definable},\ x \in \spc U}.
\]
\end{definizione}
Note that $f(x)$ is indeed a type, since for every $x \in X$, either $\phi\Pa{f(a)} \in x(a)$, or $\neg \phi\Pa{f(a)} \in x(a)$.
Since $\spc{\ }$ preserves the composition of maps, we can view $\spc{\ }$ as a covariant functor between the category of definable sets with definable maps, and the category of sets.
\begin{remark}\label{REM:IMAGE}
Let $U \subseteq A$ and $V \subseteq B$ be definable.
Then
\[
\spf^{-1}(\spc V) = \spc{f^{-1}(V)}, \quad \text{and} \quad
\spf(\spc U) = \spc{f(U)}.
\]
\end{remark}
\begin{remark}
$\spf$ is continuous with the Stone topology.
\end{remark}
\begin{remark}
If $Z \subseteq Y$  is \tdef, then $\spf^{-1}(Z)$ is quasi-compact.
\end{remark}
\begin{proof}
Since $\spf^{-1}(Z)$ is \tdef, it is also quasi-compact.
\end{proof}
\begin{remark}
If $f$ is continuous,
then $\spf$ is also continuous.
\end{remark}
\begin{proof}
Let $U \subseteq B$ be definable and open.
It suffices to prove that $\spf ^{-1}(\spc U)$ is open in~$A$.
Let $V := f^{-1}(U)$.
Since $f$ is definable and continuous, $V$~is definable and open.
Since $\spf^{-1}(\spc U) = \spc V$, we are done.
\end{proof}
Therefore, we can also view $\spc{\ }$ as a covariant functor between the category of definable sets with definable continuous maps, and the category of topological spaces.
\begin{remark}\label{REM:CONTINUOUS}
Let $x, y \in X$ such that $y \leq x$.
If $f$ is continuous, then $\spf(y) \leq \spf(x)$.
\end{remark}
\begin{proof}
Because $\spf$ is continuous.
\end{proof}
\begin{remark}
For every $x \in X$, $\dim(\spf(x)) \leq \dim x$.
\end{remark}
\begin{definizione}[Rudin-Keisler ordering]
Let $x \in \Sk(M)$ and
$y \in \Sh (M)$.
We will say that $y$ is 
less or equal to $x$ in for the Rudin-Keisler ordering, and write
$x \RKl y$, 
iff $y = \spc g(x)$ for some $g: M^k \to M^h$ definable.
We will say that $x$ and $y$ are RK-equivalent, and write $x \RKeq y$,
iff $x \RKleq y$ and $y \RKleq x$.
\end{definizione}
\begin{lemma}
Let $x \in \Types k(M)$ and $y \in \Types h (M)$.
The following are equivalent:
\begin{enumerate}
\item $y \RKleq x$;
\item for every $N \extendeq M$ elementary extension, if $N$ realizes~$x$, then $N$ realizes~$y$;
\item $M(x)$ realizes~$y$.
\end{enumerate}
\end{lemma}
\begin{proof}
$M$ has definable Skolem functions.
\end{proof}
Note that $\RKleq$ is a quasi-order, and that the relation $\RKeq$ is an equivalence relation.
\begin{lemma}\label{LEM:0-DIM}
Let $x \in \Types k(M)$ and
$y \in \Types h (M)$ be such that $y \RKleq x$.
Then, $\dim y \leq \dim x$.
Moreover, if $\dim y = \dim x$, then $y \RKeq x$.
\end{lemma}
\begin{proof}
The fact that $\dim y \leq \dim x$ is trivial.
For the second part, let $f: M^k \to M^h$ be definable such that $\spf(x) = y$.
Let $N \extendeq M$ and $c \in N^k \models x$.
Therefore, $d := f(c) \models y$, and $d \in M(c)$.
Since
\[
\dim(c/M) = \dim x = \dim y = \dim(d/M),
\]
$\dim\Pa{c/M(d)} = 0$, namely $M(c) = M(d)$.
\end{proof}
\begin{lemma}\label{LEM:BIJECTION}
Let $x \in \Sk(M)$ and $y \in \Sh(M)$ such that $x \RKeq y$, and $f: M^k \to M^h$ be definable, such that $\spf(x) = y$.
Then, there exist $A \subseteq M^k$, $B \subseteq M^h$, and $g: B \to A$ definable such that
\begin{enumerate}
\item $x \in \spc A$ and $y \in \spc B$;
\item $\dim A = \dim x = \dim y = \dim B$;
\item $f$ is continuous on $A$, and $g$ is continuous on~$B$;
\item $\spc g(y) = x$;
\item $g \comp f$ is the identity on $B$;
\item $f \comp g$ is the identity on $A$.
\end{enumerate}
\end{lemma}
\begin{proof}
Let $A \subseteq M^k$ be a cell, such that $x \in \spc A$ and $\dim A = \dim x =: n$.
\Wlog, we can assume that $n = k$ and $A = M^n$.
By Lemma~\ref{LEM:INJECTIVE}, thee exists a cell $C \subseteq M^n$ such that $x \in \spc C$ and either $f \rest C$ is injective, or $\dim f(C) < n$.
If $f \rest C$ were not injective, then $\dim y < n$,
and hence $x \RKneq y$, absurd.
Hence, $f \rest C$ is injective.

Substitute $A$ with~$C$, and let $B := f(A)$.
Hence, $f: A \to B$ is a bijection.
Let $g := f^{-1}$; the remainder of the conclusion follows (after restricting $A$and $B$ if necessary, in order to get the continuity of~$g$).
\end{proof}
\begin{lemma}\label{LEM:RETRACT-TYPE}
Let $y \in \Types k(M)$.
There exists $V \subseteq M^k$ open and definable, and a definable continuous map $\rho: V \to M^k$ such that, for every $x \in \spc V$, if $x \geq y$, then $\spr(x) = y$.
\end{lemma}
\begin{proof}
Let $E \subseteq M^n$ be a definable cell, such that $\dim E = \dim y =: n$ and $y \in \spc E$.
By Lemma~\ref{LEM:RETRACT}, there exists $V$ definable open neighborhood of $E$ and a continuous definable retraction $\rho: V \to E$.
Note that $V$ open and $y \in \spc E$ imply that $x \in \spc E$.
\begin{claim}
For every $x \geq y$, $\spr(x) = y$.
\end{claim}
By Remark~\ref{REM:CONTINUOUS}
we have $y = \spr(y) \leq \spr(x)$.
Moreover, $\spr(x) \in \spc E$.
Since $y$ is maximal in~$\spc E$, the claim is true.
\end{proof}
\begin{thm}\label{THM:REALISATION}
Let $x, y \in \Types k(M)$, such that $y \leq x$.
Then, $y \RKleq x$.
\end{thm}
\begin{proof}
Let $N \extendeq M$ be an elementary extension realizing~$x$.

Let $\rho: V \to M^k$ be as in Lemma~\ref{LEM:RETRACT-TYPE}.
Note that if $N^k \ni c \models x$, then $\spr(c) \models y$, and therefore $N$ realizes~$y$.
\end{proof}
\subsection{Closed maps}
\begin{remark}\label{REM:NORMAL}
Let $C, D \subseteq X$ be closed and disjoint.
Then, there exist $C', D' \subseteq X$ \emph{definable}, closed and disjoint such that $C \subseteq C'$ and $D \subseteq D'$.
\end{remark}
\begin{proof}
By compactness.
\end{proof}
\begin{definizione}
$f: A \to B$ is \intro{definably closed} iff for every $C \subseteq A$ closed and definable, $f(C)$ is also closed.
\end{definizione}
\begin{example}
Note that a ``definably open'' map is nothing else than an open map.
On the other hand, a (definable) map can be definably closed, without being closed.
For instance, let $M$ be countable, $A := [0,1]^2$, $B := [0,1]$, $f: A \to B$ be the projection onto the first coordinate.
Since $A$ is \dcompact, $f$~is definably closed.
Let $\eta \in Y$ be a gap (for instance, if $M = \Raz$, we can take $\eta = \frac{\sqrt{2}}{2}$).
The intervals $\lOpenInter{\eta}{1}$ and $\lOpenInter{0}{1}$ are order-isomorphic, because they are both countable, dense, and with no minimum and a maximum.
Let $g: \lOpenInter{0}{1} \to \lOpenInter{\eta}{1}$ be an order-isomorphism.
Define
\[
C := \set{(b,g(b)): 0 < b \leq 1}.
\]
Then, $C$~is closed (because $\lim_{b \to 0^+}g(b) = \eta \notin A$), but $f(C) = \lOpenInter{0}{1}$ is not closed.
\end{example}
\begin{lemma}\label{LEM:MAPS}
$f$ is an open map iff $\spc f$ is open.
$f$ is definably closed iff $\spc f$ is closed.
\end{lemma}
\begin{proof}
The ``only if'' directions are trivial.

To prove that $\spc f$ is open, it suffices to prove that for every $U \subseteq A$ open and definable, $\spc f(\spc U)$ is also open.
However, this is immediate from Remark~\ref{REM:IMAGE}.

Also immediate form the remark is the fact that if $C \subseteq A$ is closed and definable, then $\spc f(\spc C)$ is also closed (if $f$ is closed).

It remains to prove that if $C \subseteq X$ is \emph{any} closed set, then $\spc f(C)$ is closed.
Since $C$ is closed, it is \tdef, and therefore $f(C)$ is also \tdef.
Thus,
\[\begin{aligned}
C    &= \bigcap \set{\spc D: D \subseteq A \text{ definable and closed}, C \subseteq \spc D},\\
\spf(C) &= \bigcap \set{\spc E: E \subseteq A \text{ definable}, \spf(C) \subseteq \spc E}.
\end{aligned}\]
Let $\spf(C) \subseteq E$, with $E$ definable.
Then, $C \subseteq \spf^{-1}(E)$.
By compactness, there exists $D$ closed and definable such that
\[
C \subseteq D \subseteq \spf^{-1}(E).
\]
Hence, $\spf(C) \subseteq \spf(D) \subseteq E$.
However, by what we said above, $\spf(D)$ is closed.
Therefore,
\[
\spf(C) = \bigcap \set{\spf(D): D \subseteq A \text{ definable}, C \subseteq \spc E},
\]
and in particular $\spf(C)$ is intersection of closed sets, and hence closed.

\end{proof}
\begin{corollary}\label{COR:CLOSED-COMPACT}
If $A$ is \dcompact, $f$~is continuous, and $x \in X$ is a closed point, then $\spc f(x)$ is also a closed point.
\end{corollary}
\begin{proof}
The hypothesis imply that $f$ is definably closed.
\end{proof}
\begin{lemma}\label{LEM:CLOSED-POINT}
Assume that $f$ is continuous.
$f$~is definably closed iff $\spf(x)$ is a closed point for every closed point $x$.
\end{lemma}
\begin{proof}
The ``only if'' direction follows form Lemma~\ref{LEM:MAPS}.

For the other direction, let $C \subseteq A$ be definable and closed, and $D := f(C)$. We have to prove that $D$ is closed.

If, for contradiction, $D$ is not closed, let $b \in \cl D \setminus D$, and $z := \iota(b)$.
By Lemma~\ref{LEM:CLOSURE}, there exists $y \in \spc D$ such that $z < y$.

Let $y_0$ be minimal (\wrt the ordering~$\leq$) such that:
\[\begin{aligned}
y_0 \in D&\\
\exists b' \in \cl D \setminus D&\quad \iota(b') < y.
\end{aligned}\]
Choose $b_0 \in \fr D$ such that $z_0 := \iota(b_0) < y$.
Let $x \in \spf^{-1}(y_0) \cap C$ be minimal.
If $x$ is closed, $\spf(x) = y_0$ is also closed, and therefore $z_0 = y_0$, absurd.
If $x$ is not closed, then $\cl {x} = \set{x_0, \dotsc, x_n}$, where $x_0< x_1 < \dots < x_n = x \in C$, and $n \geq 1$.
Since $\spf$ is continuous, $\spf(x_i) \in \cl{\spf(x)} = \cl y$ for every $i \leq n$.
Since $\spf(x_i) \in \cl{y_0} \cap f(C)$, $f(x_i) = y_0$ for $i = 0, \dotsc, n$ by minimality of~$y_0$.
However, this contradicts the minimality of~$x$, since $x_0 < x$.
\end{proof}
\begin{example}
If $f$ is not closed, we cannot conclude that $\spf(x)$ is closed for every closed points, even if $f$ is continuous.
Here are two examples.
\begin{enumerate}
\item Let $A := \inter{0}{1}$, $B := \closedInter{0}{1}$, $f : A \to B$ be the inclusion map, and $x = 0^+$.
Then, $\spf(x) = x$ is \emph{not} closed in~$Y$, because $\clt Y{x} = \set{0,0^+}$.
However, $x$ \emph{is} closed in~$X$.
\item Let $M$ expand a real closed field, $A := M^2$, $B := M$, $\pi_i: A \to B$ $i=1,2$ be the projections on the first and second coordinate respectively.
Let $x(a_1,a_2) \in X$ be only type satisfying the following conditions:
\begin{align*}
\spc{\pi_1}(x) &= 0^+,\\
\spc{\pi_2}(x) &= 0^+,\\
a_1 \cdot a_2  &= 1.
\end{align*}
$x$~is on the infinite  branch of the hyperbola ``near infinity'', and it is closed.
However, since $\spf(x) = 0^+$, $\spf(x)$ is \emph{not} closed.
\end{enumerate}
\end{example}
\begin{example}
Assume that $f$ is not continuous.
By Lemma~\ref{LEM:MAPS}, if $f$ is definably closed, then $\spf(x)$ is closed for every closed point $x \in X$.
However, the converse is not true.
For example, let $A =\closedInter{-1}{1}$, $B = \closedInter{-1}{2}$, and $f: A \to B$ so defined:
\[
f(x) = \begin{cases}
  x     &\text{ if $x \leq 0$}\\
  x + 1 &\text{ if $x > 0$.}\\
       \end{cases}
\]
Then, $f(A) = \closedInter{-1}{0} \sqcup \lOpenInter{1}{2}$, and therefore $f$ is \emph{not} definably closed.
However, $\spf(x)$ is closed for every closed point $x \in X$.
\end{example}
\begin{corollary}\label{COR:CLOSED}
Assume that $f$ is continuous.
The following are equivalent:
\begin{enumerate}
\item $f$~is definably closed;
\item $\spf$ is closed;
\item $f(x)$ is closed for every closed point~$x$.
\end{enumerate}
\end{corollary}
\begin{proof}
$1 \Rightarrow 2$ by Lemma~\ref{LEM:MAPS}.
$2 \Rightarrow 3$ by definition.
$3 \Rightarrow 1$ by Lemma~\ref{LEM:CLOSED-POINT}.
\end{proof}
\begin{lemma}
Assume that $M$ expands a real closed field, and let $A \subseteq M^k$ be definable.
The following are equivalent:
\begin{enumerate}
\item $A\subseteq M^k$~is \dcompact;
\item for every $B$ and every $f:A \to B$ definable and continuous, $f$~is definably closed;
\item for every $B$ and every $f:A \to B$ definable and continuous, $f(A)$~is closed in~$B$.
\end{enumerate}
\end{lemma}
\begin{proof}
$1 \Rightarrow 2$ and $2 \Rightarrow 3$ are trivial.
It remains to prove that $3 \Rightarrow 1$.
Let $g: M^k \to M^{k'}$ be any definable injective continuous map, such that the image of $g$ is \emph{bounded}.
Let $B := M^{k'}$, and $f$ be the restriction of $g$ to~$A$.
Since $f(A)$ is closed and bounded, $f(A)$~is \dcompact.
Since $f$ is invertible, $A$~is also \dcompact.
\end{proof}
In the proof of the above lemma, we used the fact that $M$ expands a field only to construct the map~$g$.
However, the existence of such a map is equivalent to the fact that $M$ expands a real closed field~\cite[Corollary~9.2]{PETERZIL-STARCHENKO:1998}.


%% file: compact.tex
\section{Compactification}
\label{SEC:COMPACT}
\begin{definizione}
A \intro{\dcompactification} of $A$ (also called completion in~\cite{DRIES:1998}) is a map $\rho: A \to C$ such that:
\begin{enumerate}
\item $C \subseteq M^h$ is \dcompact;
\item $\rho$ is a definable homeomorphism onto its image;
\item $\rho(A)$ is dense in $C$.
\end{enumerate}
If the map $\rho$ is clear from the context, we will simply say that $C$ is a \dcompactification of~$A$.
\end{definizione}
\begin{example}
In the definition of \dcompactification, we cannot weaken (2) to ``$\rho$ is definable, continuous and injective''.
For instance, let
\[
A := [0,1)\subset \Real,\quad
C := \sphere^1 \subset \Complex,\quad
\rho(t) := e^{i\pi t}.
\]
$\rho: A \to C$ is not a \dcompactification.
\end{example}
\begin{lemma}
Let $\rho: A \to C$ be a \dcompactification of~$A$.
Then, $A$ is locally \dcompact iff $\rho(A)$ is open in~$C$.
\end{lemma}
\begin{proof}
For the ``only if'' direction, use Lemma~\ref{LEM:DENSE-COMPACT}.
For the other direction, let $a \in A$, and $U \subseteq C$ definable and open such that $\rho(a) \in U$ and $\cl U \subseteq \rho(A)$ ($U$~exists because $\rho(A)$ is open and $C$ is normal).
Then, $\rho^{-1}(U)$ is a relatively \dcompact neighborhood of~$a$.
\end{proof}
\begin{definizione}
Let $f : A \to B$ be continuous.
A \dcompactification $\rho: A \to C$ is \intro{compatible} with $f$ iff there exists a definable continuous map $g: C \to B$ such that $g \comp \rho = f$.
\end{definizione}
From now on, we will assume that \textbf{$M$~expands a real closed field}.
\begin{lemma}
Given a definable continuous map $f: A \to B$, where $B$ is \dcompact, there exists a \dcompactification of $A$ compatible with~$f$.%
\footnote{Thanks to M.~Mamino for the proof.
It is the same proof as in~\cite{DRIES:1998}.}
\end{lemma}
\begin{proof}
Assume that $A \subseteq M^k$ is bounded.
Let $\Gamma(f) \subseteq A \times B$ be the graph of~$f$, and
$C := \cl{\Gamma(f)}\subseteq \cl A \times B$ be its closure.
Define
\[
\fdef{\rho}{A}{C}{a}{\Pa{a,f(a)},}
\]
and $g: C \to B$ be the projection on the second coordinate.
\end{proof}
\begin{example}
Let $A := (\closedInter{0}{1} \times \closedInter{-1}{1}) \setminus{(0,0)}$, and $B := \closedInter{0}{1}$.
Define $f: A \to B$ as
\[
f(a_1, a_2) = \min \set{1, \Abs{\frac{a_2}{a_1}} }.
\]
The \dcompactification of $A$ given in the proof of the lemma is given by the disjoint union of the graph of $f$ (which is homeomorphic to~$A$), and the vertical segment $\set{(a_1,a_2,b) \in M^3: a_1 = a_2 = 0, 0 \leq a_3 \leq 1}$.
\end{example}
\begin{remark}
If $\dim A = 1$, there exists a \intro{universal \dcompactification} of~$A$ (namely, one compatible with all the definable functions $f$ with \dcompact co-domains).
\end{remark}
\begin{example}
If $\dim A > 1$, such universal \dcompactification might not exist.
Let $A := \lOpenInter{0}{1} \times \closedInter{-1}{1}$.
A universal \dcompactification for $A$ does not exist.
\end{example}
\begin{definizione}
Let $x \in X$ and $C \subseteq X$ be definable.
We shall say that $x$ is \intro{near}~$C$, and write $C \leq x$, iff every definable open neighborhood of $C$ contains~$x$.
We shall write $C < x$ iff $C \leq x$ and $x \notin \spc C$.
\end{definizione}
\begin{lemma}\label{LEM:NEAR}
The following are equivalent:
\begin{enumerate}
\item $C \leq x$;
\item every \emph{definable} closed set containing $x$ intersects~$C$;
\item every closed set containing $x$ intersects~$C$;
\item there exists $y \in C$ such that $y \leq x$ (namely, $\cl x \cap C \neq \emptyset$).
\end{enumerate}
\end{lemma}
\begin{proof}
$1 \Leftrightarrow 2$ and $3 \Rightarrow 2$ are trivial.

For $2 \Rightarrow 4$, let
\[
\Cfam := \set{D \cap C: D \subseteq X \text{ closed and definable},\
x \in D}.
\]
$\Cfam$ is a family of definable subsets of $C$ with the \FIP.
Since $C$ is compact (with the Stone topology),
$\bigcap \Cfam \neq \emptyset$.
Any $y \in \bigcap \Cfam$ is in $\cl x \cap C$.


For $4 \Rightarrow 3$, any closed set $D$ containing $x$ must also contain $y$, and therefore \mbox{$C \cap D \neq \emptyset$}.
%
\end{proof}
\begin{remark}
\[
\set{x \in X: C \leq x} = \bigcap \set{U: C \subseteq U \subseteq X,
U \text{ open and definable}}.
\]
\end{remark}
The above lemma suggests the following extension of the notion $C \leq x$ to the case when $C$ is \tdef.
\begin{lemma}
Let $Z \subseteq X$ be \tdef. Then,
\[\begin{aligned}
\hat Z &:= \bigcap \set{U: Z \subseteq U \subseteq X,\
U \text{ open and definable}}=\\
&=\bigcap \set{U: Z \subseteq U \subseteq X,\
U \text{ open}}=
\set{x \in X: \cl x \cap Z \neq \emptyset}.
\end{aligned}\]
We say that $x$ is near $Z$, and write $Z \leq x$, iff $x \in \hat Z$.
\end{lemma}
Note that $\hat Z$ is \tdef.
\begin{definizione}
Let $x \in X$.
We will say that $x$ is far from the frontier of $X$ iff there exists a definable set $C \subseteq A$ which is \dcompact and such that $x \in \spc C$.
Otherwise, we will say that $x$ is near the frontier.
The \intro{fringe} $\frin A$ of $A$ is the set of points of $X$ near the frontier.
\end{definizione}
\begin{lemma}
The following are equivalent:
\begin{enumerate}
\item $x \in \frin A$;
\item for \emph{every} $C$ \dcompactification of~$A$, $\fr A < x$, where $\fr A$ is the frontier of $A$ taken inside~$C$;
\item for \emph{some} $C$ \dcompactification of~$A$, $\fr A < x$;
\item for every definable $D \subseteq A$, if $D$ is \dcompact, then $x \notin \spc D$.
\end{enumerate}
\end{lemma}
\begin{proof}
Assume that $x\in \frin A$.
Let $U$ be an open definable neighborhood of $\fr A$,
and $D := C \setminus U$.
Since $C$ is \dcompact, $D$~is also \dcompact and contained in~$A$, and therefore $x \notin \spc D$.
Hence, $x \in \spc U$, and thus $\fr A \leq x$.
Since $x \in X$, $x \notin \fr A$, and so \mbox{$\fr A < x$}.

Conversely, assume that $\fr A < x$ for some $C$ \dcompactification of~$A$.
Let $D \subseteq A$ be \dcompact, and $U := A \setminus D$.
Since $D$ is \dcompact, $U$~is open in $A$.
$D \cap \fr A = \emptyset$, because $D \subseteq A$, and so $U$ is an open neighborhood of $\fr A$.
Thus, $x \in \spc U$, and therefore \mbox{$x \in \frin A$}.
\end{proof}
\begin{definizione}\label{DEF:SCLOSED}
Let $x \in X$.
We will say that $x$ is \intro{strongly closed} iff $x$ is closed and $x$ is far from the frontier of~$A$.
\end{definizione}
\begin{lemma}
The following are equivalent:
\begin{enumerate}
\item $x$ is strongly closed;
\item $x$ is closed in some \dcompactification of~$A$;
\item $x$ is closed in every \dcompactification of~$A$;
\item $x \in \spc D$, for some $D \subseteq A$ definable and \dcompact.
\end{enumerate}
\end{lemma}
\begin{proof}
$1 \Rightarrow 3 \Rightarrow 2$ is trivial.
For $2 \Rightarrow 1$, assume that $x$ is closed in some \dcompactification $C$ of $A$. It is then trivially true that $x$ is already closed in~$A$.
If, for contradiction, $\fr A < x$, we would have that there exists $y \in \fr A$ such that $y < x$, contradicting the fact that $x$ is closed in~$\spc C$.
\end{proof}
\begin{remark}
If $A$ is \dcompact, then $x$ is \sclosed iff it is closed.
\end{remark}
\begin{lemma}\label{LEM:SCLOSED}
Let $A \subseteq B$, and $x \in X$.
Then, $x$ is strongly closed in $A$ iff it is strongly closed in~$B$.
\end{lemma}
\begin{proof}
Let $\rho: B \to D$ be a \dcompactification of $B$.
\Wlog, we can assume that $\rho$ is the inclusion map.
Let $C := \clt D A$.
Then, $C$ is a \dcompactification of~$A$.

If $x$ is not \sclosed in~$A$, then there exists $y \in \spc C$ such that $y < x$. Since $C \subseteq D$, $x$ is not closed in $D$ either, and therefore $x$ is not \sclosed in $B$.

If $x$ is not \sclosed in $B$, then there exists $y \in \spc D$ such that $y < x$.
However, $C$ is closed in $D$, and therefore $y \in \spc C$.
Hence, $x$ is not \sclosed in $A$ either.
\end{proof}
\begin{thm}\label{THM:SCLOSED-MAP}
Let $f : A \to B$ be \emph{any} definable map.
If $x \in X$ is strongly closed (in~$A$), then $\spf(x)$ is also strongly closed (in~$B$).
\end{thm}
\begin{proof}
Decompose $A$ into cells $C_i$ such that $f$ is continuous on each cell.
Let $C$ be the cell counting~$x$.
By Lemma~\ref{LEM:SCLOSED}, $x$ is strongly closed in~$C$.
Let $E$ be some \dcompactification of $B$.
Consider the map $g: C \to E$ given by the composition of the immersion of $B$ in $E$ with the restriction of $f$ to $C$.
Let $\rho: C \to D$ be a \dcompactification of $C$ compatible with~$g$.
By Corollary~\ref{COR:CLOSED-COMPACT}, $\spc g(x) = \spf(x)$ is closed in~$E$, and therefore $\spf(x)$ is strongly closed.
\end{proof}
In a slogan, ``being strongly closed'' is an intrinsic property of a type $x$ (namely, independent from the ambient space~$A$).

%% file: definable.tex
\section{Rational and irrational types}
\label{SEC:RATIONAL}
\begin{definizione}
A type $x\in\spc{A}$ is \intro{rational} (a.k.a.~definable) iff for every formula $\psi(a,u)$ there exists $c \in A$ such that
\[
\set{b \in M^h: \psi(a,b) \in x(a)} =
\set{b \in M^h: M \models \psi(c,b)}.
\]
If $x$ is not rational, then it is \intro{irrational}.
\end{definizione}
Note that $\iota(a)$ is a \ratt type for every $a \in A$.

Remember that $M$ expands a field.
\begin{lemma}
Let $N$ be an elementary extension of~$M$.
Let $A$ be a \dcompactification of~$M$ (\eg, $A = \sphere^1$).
The following are equivalent:
\begin{itemize}
\item $x$~is rational (and not realized);
\item either $x = \pm \infty$, or there exists $a \in M$ such that $x = a^{\pm}$;
\item $M$ is \emph{not} co-final in~$M(x)$ 
\item $M$ is Dedekind complete in $M(x)$;
\item $x$ is not closed in $\spc A$;
\item there exists $a \in A$ such that $\iota(a) < x$.
\end{itemize}
The following are also equivalent:
\begin{itemize}
\item $x$ is irrational;
\item $M$ is co-final in $M(x)$;
\item $x$ is closed in $\spc A$.
\end{itemize}
\end{lemma}
\begin{proof}
See \cite{PILLAY:1994} (or \cite{BAI-POIZAT:1998}).
Note that the hypothesis imply that $\dim x = 1$.
\end{proof}

The following example is by~Coste.
\begin{example}\label{EX:GAP-EXP}
Let $M' := \Pa{\Real, + , \cdot, \exp}$, and $A = [0,1]^2$.
Let $x = x(a_1,a_2) \in \spc A$ be the unique type satisfying the conditions
\begin{align*}
a_1 &\models 0^+\\
a_2 &= \exp (a_1).
\end{align*}
Note that $a_2 \models 0^+$.
Moreover, $\dim(x) = 1$, and $0 < x$.

Let $M$ be the reduct of $M'$ to the field structure alone, and $y$ be the image of $x$ under the reduct map.
Note that $0 < y$, and that $\dim y = 2$, because the germ of $\exp$ near $0$ is not definable in~$M$.
Moreover, there is no $z \in \spc M$ such that $0 < z < y$, otherwise there would be $z' \in \spc {M'}$ such that $0 < z' < x$, which is impossible.
Note also that $y$ is a rational type, because \emph{all} type over $M$ are rational, since $M$ is Dedekind complete.
\end{example}
\begin{example}\label{EX:GAP-E}
Let $M$ be the field of real algebraic numbers.
Let $N$ be the real closure of $M(\varepsilon, e)$, where $\varepsilon$ is a positive infinitesimal element (namely, $\tp \varepsilon M = 0^+$),
and $e$ be a transcendental real number.
Let $A := [0,1]^2 \subseteq M^2$, and $x(a_1,a_2) \in \spc A$ be given by
\begin{align*}
a_1 &\models 0^+\\
a_2 &= e a_1.
\end{align*}
Note that $a_2 \models 0^+$, that $0 < x$, that $\dim x = 2$, but there is no $y \in \spc A$ such that $0 < y < x$.
Moreover, $x$~is irrational, because the type $\tp e M$ can be defined using~$x$, and $\tp e M$ is irrational.
\end{example}
Let $M \elemeq N$ be an elementary extension.
\begin{definizione}
We will say that the extension $N/M$ is \intro{rational} iff, for every $n\in\Nat$, every $n$-type over $M$ realized in $N$ is rational.
\footnote{Rational extensions were called \emph{tame} extensions in~\cite{DRIES-LEW:1995} (in the case of o-minimal theories expanding a real closed field).}
We will say that $N/M$ is \intro{\tirrational} iff , for every $n\in\Nat$, and for every $c \in N^n \setminus M^n$, $\tp c M$ is irrational.
We will say that a type $x \in \Sn(M)$ is \tirrational iff there exists a totally irrational extension $N$ realizing~$x$, or equivalently iff $M(x)/M$ is \tirrational.
\end{definizione}
\begin{remark}
A \tirrational type is either irrational, or already realized in~$M$.
\end{remark}
\begin{lemma}
$N/M$ is rational iff every $1$-type realized in $N$ is rational.
$N/M$ is \tirrational iff, for every $c \in N \setminus M$, $\tp c M$ is irrational.
\end{lemma}
\begin{proof}
See \cite{PILLAY:1994} or \cite{BAI-POIZAT:1998}.
\end{proof}
\begin{remark}
$N/M$ is \tirrational iff $M$ is \emph{co-final} in~$N$.
\end{remark}
\begin{remark}
If $M \elemeq N \elemeq P$ and $N/M$ is \tirrational, then $\mequiv$ and $\nequiv$ coincide (on $P^k$).
\end{remark}
\begin{definizione}
Let $M \elemeq N$, and $b,c \in N^k$.
We shall say that $b$ is \intro{\mbounded}, (or simply bounded if $M$ is clear from the context) iff there exists $a \in M$ such that $\abs{b} \leq a$.
We shall write $b \mequiv c$ iff for every $a \in M$ such that $a > 0$, we have $\abs{b-c} < a$.
Note that $\mequiv$ is an equivalence relation.
\end{definizione}
\begin{lemma}
$N/M$ is rational iff for every $k \in \Nat$ and every \mbounded $c \in N^k$ there exists (a necessarily unique) $a \in M^k$ such that $a \mequiv c$.
\end{lemma}
We will call $a$ as above the \intro{\mstandard part} of~$c$, and write $a = \stm(c)$ (or simply $a = \st(c)$ if $M$ is clear from the context).
\begin{remark}\label{REM:TIRRATIONAL}
$x \in \Sn(M)$ is \tirrational iff, for every definable function \mbox{$f: M^n \to [0,1]$}, we have $\spf(x) \neq 0^+$.
\end{remark}
\begin{lemma}\label{LEM:DISTANCE}
Let $C \subseteq A$ be definable.
Define $f: A \to M$ as $f(a) := d(a,C)$, where $d$ is the distance.
Let $x \in X$ and $y := \spf(x) \in \spc{M}$.
Then, $C < x$ iff $y = 0^+$.
\end{lemma}
\begin{proof}
For every $\varepsilon \in M^{>0}$, let
$V_\varepsilon := f^{-1}\rOpenInter{0}{\varepsilon}$.

Assume that $C < x$. Since $x \notin \spc C$, $y \in \spc{\inter{0}{+\infty}}$.
Since, $\forall \varepsilon > 0$, $V_\varepsilon$ is an open definable neighborhood of~$C$, and $C \leq x$, we have $x \in \spc {V_\varepsilon}$, and therefore $y \in \spc{\inter{0}{\varepsilon}}$.
Thus, \mbox{$y \in \bigcap_{\varepsilon > 0} \spc{\inter{0}{\varepsilon}}$}, and so $y = 0^+$.

Conversely, if $y = 0^+$, let $U$ be a definable open neighborhood of~$C$.
By Lemma~\ref{LEM:NEIGH-F}, $V_\varepsilon \subseteq U$ for some $\varepsilon > 0$. Since $y = 0^+$, $x \in V_\varepsilon$.
Therefore, $C \leq x$.
Since $y \neq 0$, $x \notin \spc C$, and therefore $C < x$.
\end{proof}
\begin{thm}\label{THM:SCLOSED}
Let $x \in \Sn(M)$.
Then, $x$ is \sclosed iff it is \tirrational.
\end{thm}
\begin{proof}
Let us prove \textbf{the ``if'' direction}.
Assume that $x$ is not \tirrational.
Then, by Remark~\ref{REM:TIRRATIONAL}, there exists $f: M^n \to [0,1]$ definable such that $\spf(x) = 0^+$.
Decompose $M^n$ into cell such that $f$ is continuous on each cell, and let $C$ be the cell containing~$x$.
Let $\rho: C \to D$ be a \dcompactification of $C$ compatible with~$f$, and let $g: D \to \closedInter{0}{1}$ be the corresponding extension of~$f$.
\Wlog, $\rho$ is the inclusion map.
Moreover, $\spc g(x) = 0^+$.
Let $E := g^{-1}(0)$.
\begin{claim}
$E \leq x$.
\end{claim}
Let $U$ be a definable neighborhood of~$E$.
By Lemma~\ref{LEM:NEIGH-F}, there exists $\varepsilon > 0$ such that $g^{-1}\Pa{\inter{0}{\varepsilon}} \subseteq U$.
Since $\spc g(x) = 0^+$, we have
\[
\spc x \in g^{-1}\Pa{\inter{0}{\varepsilon}} \subseteq U,
\]
and the claim is proved.

Since $x \notin E$, we have $E < x$.
By Lemma~\ref{LEM:NEAR}, there exists $y \in E$ such that $y < x$.
Therefore $x$ is not closed in~$D$, and thus $x$ is not \sclosed.

For \textbf{the ``only if'' direction},
assume that $x$ is not \sclosed.
\Wlog, we can assume that $x \in \spc A$ for some \dcompact set~$A$.
Let $y < x$, and let $C \subseteq A$ be a closed definable \subseT such that $y \in \spc C$ and $\dim C < \dim x$.
By Lemma~\ref{LEM:NEAR}, $C \leq x$, and, since $\dim C < \dim x$, $C < x$.
Therefore, by Lemma~\ref{LEM:DISTANCE}, $\spf(x) = 0^+$, where $f(a) := d(a,C)$.
Thus, by Remark~\ref{REM:TIRRATIONAL}, $x$~is not \tirrational.
\end{proof}
\begin{lemma}\label{LEM:LIFTING}
Let $A$ be \dcompact.
Let $f: A \to B$ be definable and continuous, $x \in \spc A$ and $y \in \spc B$ such that $y \leq \spf(x)$.
Then, there exists $z \in \spc A$ such that $z \leq x$ and $\spf(z) = y$.
Moreover, if $y < \spf(x)$, then $z < x$.
\end{lemma}
\begin{proof}
Let $C := \cl x$.
Since $A$ is \dcompact, $\spf$ is closed, and therefore $\spf(C)$ is closed.
Since $x \in \spf(C)$ and $y \leq x$, we have $y \in \spf(C)$.
Let $z \in C$ such that $\spf(z) = y$.
\end{proof}
\begin{example}
We cannot drop the condition that $A$ is \dcompact in Lemma~\ref{LEM:LIFTING}.
For instance, let $B := [0,1]$, $A := B \setminus\mset{0}$, $x := 0^+$, $y := 0$ and $f: A \to B$ be the inclusion map.
A point $z$ as in the conclusion of the Lemma does not exists.
\end{example}
\begin{remark}
Let $y \RKleq x$.
If $x$ is rational, then $y$ is also rational.
If $x$ is \tirrational, then $y$ is also \tirrational.
\end{remark}
\begin{lemma}\label{LEM:STEP-1}
Let $y < x \in X$.
Assume that $x$ is rational and $n := \dim x = 1 + \dim y$.
Then, $y$~is rational.
\end{lemma}
\begin{proof}
\Wlog, we can assume that $A \subseteq [0,1]^k$.
We will prove the conclusion by induction on~$k$.
The case $k = 0$ is trivial, and the one $k = 1$ is easy.
Let $f: M^k \to [0,1]$ be definable.
We want to prove that $\spf(x)$ is rational.
Decompose $M^k$ into cells in a way compatible with $A$ and with~$f$.
Let $C \subseteq M^k$ be the cell containing~$x$.
\footnote{Note that $C$ might not contain~$y$.}

Let $\rho: C \to D$ be a \dcompactification of $C$ compatible with $f$ and with the inclusion map $\lambda: C \to [0,1]^k$, and $g: D \to [0,1]$, $\mu: D \to [0,1]^k$ be the extensions of $f$ and $\lambda$ respectively.
If we identify $C$ with $\rho(C)$, we can assume that $\rho$ is the inclusion map.
By Lemma~\ref{LEM:LIFTING}, there exists $z < x$ such that $\spc\mu(z) = y$.
\begin{claim}
$\dim z = \dim y$.
\footnote{Here is the point where we use $\dim y = n-1$.}
\end{claim}
In fact, $z < x$ implies that $\dim z \leq n-1$.
Moreover, since $y = \spc\mu(z)$, $\dim z \geq \dim y$.
\begin{claim}\label{CL:RATIONAL}
$z$ is rational.
\end{claim}
In fact, by Lemma~\ref{LEM:0-DIM}, $z \RKleq y$, and $y$ is rational.

By Claim~\ref{CL:RATIONAL}, $\spc g(z)$ is \ratt.
Moreover, since $g$ is continuous, $\spc g(z) \leq \spc g(x)$.
By the case $k = 1$, $\spc g(x)$ is also \ratt.
\end{proof}
\begin{corollary}\label{COR:CHAIN}
Let $x_0, \dotsc, x_n \in A$ such that $x_0 < x_1 < \dots < x_n$ and $\dim(x_n) = n$.
Then, each $x_i$ is rational, for $i = 0, \dotsc, n$.
\end{corollary}
\begin{proof}
Note that $\dim (x_i) = i$ for every $i \leq n$.
Since $\dim x_0 = 0$, $x_0$ is rational.
By induction on $n$, $x_{n-1}$ is rational.
By Lemma~\ref{LEM:STEP-1}, $x_n$ is also rational.
\end{proof}
\begin{lemma}\label{LEM:RATIONAL-0}
Let $A$ be \dcompact, and $x \in X$.
If $x$ is rational, then there exists a unique $a \in A$ such that $\iota(a) \leq x$.
\end{lemma}
\begin{proof}
Uniqueness is trivial.
Assume, for contradiction, that $a$ does not exist.
Let $y$ be the minimum of $\cl x$. 
Note that $y$ is closed in~$A$.
Since $A$ is \dcompact, $y$~is \sclosed, and therefore, by Theorem~\ref{THM:SCLOSED}, $y$ is \tirrational.
Since $\dim y > 0$, $y$ is irrational.
Since $y \leq x$, by Theorem~\ref{THM:REALISATION}, $y \RKleq x$, and therefore $x$ is irrational, absurd.
\end{proof}
\begin{definizione}
Let $A$ be \dcompact, and $x \in X$.
We will say that $\rho : C \to D$ is a \dcompactification fixing~$x$ iff
\begin{enumerate}
\item $C \subseteq A$ is definable, such $x \in \spc C$;
\item $\rho : C \to D$ is a \dcompactification compatible with the inclusion map $\lambda : C \to A$.
\end{enumerate}
In this case, we will denote by $\mu: D \to A$ the extension of~$\lambda$ to~$D$.
\end{definizione}
\begin{definizione}
Let $A$ be \dcompact.
Let $x,y \in \spc A$ such that $y \leq x$.
Let $\rho : C \to D$ be a \dcompactification fixing~$x$.
By Lemma~\ref{LEM:LIFTING}, there exists $z \in D$ such that $z \leq \spr(x)$ and $\spm(z) = y$.
We will call the pair $(z,\spr(x))$ a \intro{lifting} of $(y, x)$ compatible with the \dcompactification~$\rho$.
We will say that the pair $(y,x)$ is \intro{\maximal} iff for every lifting $(z,x')$ of $(y,x)$, $\dim z \leq \dim y$.
\end{definizione}
Note that if $(z,x')$ is a lifting of $(y,x)$, then $\dim y \leq \dim z$.
Therefore, if $(y,x)$ is \maximal, we have $\dim y = \dim z$, and hence, by Lemma~\ref{LEM:0-DIM}, $y \RKeq z$.
\begin{lemma}\label{LEM:UNIQUENESS}
Let $A$ be \dcompact.
Let $y \leq x \in X$, with $(y,x)$ \maximal.
Let $\rho: C \to D$ be a \dcompactification fixing~$x$, and
$x' := \spr(x)$.
Then, there exists a \emph{unique} $z \in \spc D$ such that $z \leq x'$ and $\spm(z) = y$.
\end{lemma}
\begin{proof}
The existence of $z$ is Lemma~\ref{LEM:LIFTING}.
For the uniqueness, let $z_1 \leq z_2 \leq x'$ such that $\spm(z_1) = \spm(z_2) = y$.
However, by \maximality of $(x,y)$, $\dim z_1 = \dim y = \dim z_2$.
Thus, $z_1 = z_2$.
\end{proof}
Therefore, if $(y,x)$ is \maximal, and $\rho$ is as in the hypothesis of the lemma, we can speak of \emph{the} lifting of $(y,x)$ compatible with~$\rho$.
\begin{definizione}
Let $x \in \Sk(M)$, $y \in \Sh(M)$, such that $y \RKleq x$,
and $f: M^k \to M^h$ be a definable function, such that $\spf(x) = y$.
We will say that $(y,x,f)$ is \intro{rational} iff, for every $N \extendeq M$, and every $N^k \ni c\models x$, the type $\tp c {N(d)}$ is rational, where $d := f(c)$.
\end{definizione}
\begin{remark}
Let $x, y, f$ be as in the above definition.
$(y,x,f)$ is rational iff, for \emph{some} $N \extendeq M$ and some $N^k \ni c \models x$, the type $\tp c {N(d)}$ is rational, where $d := f(c)$.
\end{remark}
\begin{remark}
Note that if $y = \iota(a)$ for some $a \in M^k$, then $(y,x,f)$ is rational
(for \emph{any} $f: M^k \to M^h$ such that $\spf(x) = y$) iff $x$ is
rational.
\end{remark}
\begin{example}
The fact that $(y,x,f)$ is rational does depend on the particular choice of~$f$ (such that $\spf(x) = y$).
For instance, let $k := 2$, $h := 1$,
\[
x(a_1,a_2) := a_i \models 0^+, a_1 \ll a_2,
\]
$y := 0^+$, $f_1(a_1,a_2) := a_1$, $f_2(a_1,a_2) := a_2$.
Then, $(y,x,f_1)$ is not rational, while $(y,x,f_2)$ is rational.
\end{example}
\begin{remark}
If $\dim x = 1$, the fact that $(y,x,f)$ is rational does \emph{not} depend on~$f$.
\end{remark}
\begin{proof}
Let $N \extendeq N$, and $N \ni c \models x$, and $d := f(c) \models y$.
If $d \in M$, then $(y,x,f)$ is rational iff $x$ is rational.
If $d \notin M$, then $M(d) = M(c)$, and therefore $(y,x,f)$ is always rational.
\end{proof}
The following Lemma is contained in~\cite[Proposition~2.1]{PILLAY:1988}.
\begin{lemma}\label{LEM:EXTENSION}
Let $N \extendeq M$, $d \in N^k$, $N' := M(d)$, $e \in N'^h$.
Then, there exist $U$ neighborhood of $d$ definable (in~$M$), $h: U \to N^h$ definable and continuous, such that $h(d) = e$.
\end{lemma}
\begin{proof}
Since $e \in N'^h$, there exists $h': N^k \to N^h$ definable (but not necessarily continuous) such that $h'(d) = e$.
Let $C \subseteq N^k$ be a cell (definable in $M$!) such that $h'$ is continuous on $C$ and $d \in C$.
By Lemma~\ref{LEM:RETRACT}, there exists an open neighborhood $U$ of $C$ and retraction $\rho: U \to C$.
Define $h := h' \comp \rho: U \to N^h$.
\end{proof}
\begin{thm}\label{THM:MAX-RAT}
Let $x, y \in \Sk(M)$, with $y \leq x$.
By Lemma~\ref{LEM:RETRACT-TYPE},
there exists
$f: U \to C$ retraction such that $\spf(x) = y$, where $C$ contains~$y$, and $U$ is a neighborhood of~$C$.
If $(y,x)$ is \maximal, then $(y,x,f)$ is rational.
\end{thm}
\begin{proof}
\Wlog, we can assume that $C \subseteq A$, where $A = [0,1]^k$.
Let $N \extendeq M$ and $c \in A(N)$ such that $c \models x$, and let $d := f(c) \models y$.
Let $N' := M(d)$.
Assume, for contradiction, that there exists $g: A(N) \to N$ definable in $N'$, such that $z := \tp{g(c)}{N'}$ is not rational.
Note that $g(t) = g'(t,d)$  for some $g' : A \to M$ definable in~$M$.
Hence, $g(c) = g'(c,f(c))$, and therefore \wloG we can assume that $g$ is definable in~$M$.
Moreover, \wloG we can assume that $g(d) < g(c)$.

Let $E \subseteq U$ be a cell such that $g$ is continuous on $E$ and $x \in \spc E$.
Let $\rho: E \to D$ be a \dcompactification of $E$ compatible with $f$, $g$, and the embedding $\lambda: E \to A$.
Let $(y',x')$ be the lifting of $(y,x)$ compatible with~$\rho$.
By hypothesis, $\dim y' = \dim y$, and therefore $y \RKeq y'$.
Also, $x \RKeq x'$.

Let $\mu: D \to A$ be the extension of~$\lambda$.
By Lemma~\ref{LEM:BIJECTION}, there exist $F \subseteq E$, $G \subseteq D$ and $\nu: G \to A$ continuous, all definable, such that $y \in \spc F$, $y' \in \spc G$ and $\nu$ is the inverse of $\mu\rest F$.
Define $f' := \nu \comp f \comp \mu$; the domain of $f'$ is $U':= \mu^{-1}\Pa{f^{-1}(G)} \subseteq U$: note that $x', y' \in \spc U'$; its co-domain is~$F$.
Note that $f'$ is continuous.
Moreover,
\[
f' \comp f' = \nu \comp f \comp \mu \comp \nu \comp f \comp \mu = \nu \comp f \comp \mu = f'.
\]
Let $c' = \rho (c) \models x'$ and $d' := f'(d') \models y'$.
It suffices to prove that $(y',x', f')$ is rational.
Therefore, \wloG we can assume that $y \in \spc E$ (and $x = x'$, $y = y'$, $f =f'$, $d = d'$, $c = c'$).


Since $z$ is not rational, there exists $h : E \to M$ definable, such that
\[
g(d) < h(d) < g(c).
\]
By Lemma~\ref{LEM:EXTENSION}, we can assume that $h$ is defined and continuous on all~$E$.
Let
\[
F := \set{a \in U \cap E: h(f(a)) \leq g(a)}
\]
Note that $c \in  F(N)$, and therefore $x \in \spc F$.
Moreover, $h,f,g$ are continuous, hence $F$ is closed in~$E$.
Thus, $y \in \spc F$, and therefore
\[
h(f(d)) \leq g(d).
\]
Since $f(d) = d$, we have $h(d) \leq g(d)$, a contradiction.
\end{proof}
\begin{corollary}\label{COR:MAX-RAT}
Let $a \in A$ and $x \in X$ such that $\iota a \leq x$.
If $(\iota a, x)$ is maximal, then $x$ is rational.
\end{corollary}
\begin{remark}
Let $y \leq x \in X$.
Then, there exists a lifting $(y',x')$ of $(y,x)$, such that $(y',x')$ is maximal.
\end{remark}
\begin{paragraph}{Lifting the closure of a type.}
Let $A$ be \dcompact, and $x \in X$, with $\dim x = n$.
We shall write $\cl{x} = (x_0 < x_1 < \dotsb < x_m)$ iff
$\cl{x} =: \set{x_0, x_1, \dotsc, x_m}$, with $x_0 < x_1 \dotsb < x_m = x$.

Let $\cl{x} = (x_0 < \dotsb < x_m)$.
Let $\rho_i: C_i \to D_i$ be \dcompactifications fixing $x$, $i = 0,1, 2$, such that
\begin{enumerate}
\item $C_0 \subseteq C_1 \cap C_2$;
\item $\rho_0$ is compatible with $\rho_1$ and $\rho_2$;
\end{enumerate}
($\rho_0$ is a common refinement of $\rho_1$ and~$\rho_2$).
Let $\mu_i: D_i \to A$ be the extension of the inclusion $\lambda_i: C_i \to A$, $i = 0,1, 2$, and $\nu_i: D_0 \to D_i$ be the extension of~$\rho_i$, $i=1, 2$.
For $i = 0,1,2$, let $y_i := \spc{\rho_i}(x) \in \spc{D_i}$, and
$\cl{y_i} =: \set{z_{i,0}, \dotsc, z_{i,m_i}}$, with $z_{i,0} < z_{i,1} , \dotsb < z_{i,m_i} = y_i$.
Note that, for each $i \leq 3$ and $j\leq m_i$,
\[
\exists! f(i,j) \leq m \quad \spc{\mu_i}(z_{i,j}) = x_{f(i,j)}.
\]
Moreover, $\dim x_{f(i,j)} \leq \dim z_{i,j}$.
Similarly, for $i=1,2$ and $j \leq m_0$,
\[
\exists! g(i,j) \leq m_i \quad \spc{\nu_j}(z_{0,j}) = z_{i,g(i,j)}, \quad \text{and} \quad
\dim z_{i,g(i,j)} \leq \dim z_{0,j}.
\]

We will say that $\cl{x}$ is \intro{\cmaximal} iff for every \dcompactification $\rho: C \to D$ fixing~$x$, if we call $x' := \spr(x')$ and $\cl{x'} = (x'_0 < \dotsb < x'_{m'})$, then $m = m'$, $\spr(x'_i) = x_i$, and $\dim(x'_i) = \dim(x_i)$, for $i = 0, \dotsc, m$.

Therefore, there exist a \dcompactification $\rho: C \to D$ fixing~$x$, such that $\cl{\spr(x)}$ is \cmaximal.

Note that $x_0$ is always \sclosed (because $A$ is \dcompact).
If moreover $\cl{x}$ is \cmaximal, then $(x_i,x_j)$ is \maximal for every $i \leq j \leq m$, and in particular $(x_0,x)$ is \maximal.
\end{paragraph}
\begin{corollary}
Let $N = M(c)$, for some finite tuple~$c \in [0,1]^k(N)$.
Let $x := \tp c M$.
Let $y \leq x$ such that $y$ is \sclosed and $(y,x)$ is maximal (by the above discussion, we can always assume that $y$ exists).
Let $N \ni d \models y$, and $N' := M(d)$.
Then, $N'/M$ is \tirrational, and $N/N'$ is rational.
\end{corollary}

%% file: amalgam.tex
\section{Amalgamation}
\label{SEC:AMALGAM}
\begin{lemma}
Let $M \elemeq N$.
There exists $N'$ such that $M \elemeq N' \elemeq N$, $N'/M$ is \tirrational an $N/N'$ is rational.
There exists $N''$ such that $M \elemeq N'' \elemeq N$, $N''/M$ is rational an $N/N''$ is \tirrational.
\end{lemma}
\begin{proof}
See \cite[Lemma~3]{BAI-POIZAT:1998}.
\end{proof}
\begin{lemma}
Let $M \elemeq N \elemeq P$.
\begin{enumerate}
\item ($N/M$ is \tirrational and $P/N$ is \tirrational) iff $P/M$ is \tirrational.
\item If $N/M$ is rational and $P/N$ is rational, then $P/M$ is rational.
\item If $P/M$ is rational (resp. \tirrational), then $N/M$ is rational (resp. \tirrational).
\end{enumerate}
\end{lemma}
\begin{example}
It is not true that $P/M$ rational implies $P/N$ rational.
For instance, let $P := M(b,c)$, where $\tp b M = 0^+$
and $\tp c {M(b)} = 0^+$.
Let $N := M(c)$.
Then, $P/N$ is not rational (in fact, it is \tirrational).

Another example is in~\cite{BAI-POIZAT:1998}, where $P$ is ``very saturated'' over~$M$.
\end{example}
\begin{lemma}\label{LEM:EXT-TYPE}
Let $N \extendeq M$.
Let $N_i$, $i = 1,2$, such that
\begin{enumerate}
\item $M \elemeq N_i \elemeq N$;
\item $N_i /M$ is \tirrational;
\item $N / N_i$ is rational.
\end{enumerate}
Then, $N_1$ and $N_2$ satisfy the same types over~$M$.
\end{lemma}
\begin{proof}
Let $d \in N_1^k$, and $x := \tp{d}{M}$.
We have to prove that $x$ is realized in~$N_2$.
\Wlog, we can assume that $d \in [0,1]^k(N)$.
By Theorem~\ref{THM:SCLOSED}, $x$ is \sclosed.

Let $y := \tp{d}{N_2}$, and $z \leq y$ minimum.
Note that $z$ is \sclosed, and that, since $d\models y$, we have that $z$ is realized in~$N$.
Since, by hypothesis, all types over $N_2$ realized in $N$ are rational, $z$ is both \tirrational and rational, and therefore $z$ is realized in $N_2$.
\end{proof}
\begin{proposition}\label{CONJ:ISOM}
If $M, N, N_i$ are as in the above lemma, then $N_1$ and $N_2$ are isomorphic over~$M$.
\end{proposition}
\begin{proof}
It was proved in~\cite[Theorem~2.15]{DRIES-LEW:1995}.
The proof goes as follows.
\begin{claim}
For every $b \in N_1$ there exists a unique $c \in N_2$ such that
$b \mequiv c$.
\end{claim}
The uniqueness follows form the fact that $M$ is cofinal in $N_2$.
For the existence, if $b$ did not exist, then $\tp c {N_2}$ would not be rational, and hence $N_2(c)/N_2$ would be \tirrational, contradicting the maximality of $N_2$.
\begin{claim}
The map sending $b \in N_1$ to the above $c \in N_2$ is an $M$-isomorphism.
\end{claim}
Working by induction, we can assume that $N_1 = M(b)$ and $N_2 = M(c)$.
If $b \in M$, then $c = b$, and we are done.
If $b \notin M$, then $c$ fills the same gap on $M$ as $b$, and therefore $M(b)$ is isomorphic to $M(c)$ over~$M$.
\end{proof}
\begin{definizione}
Let $x \in \spc{A}$ and $y \in \spc{B}$.
Define
\[
x \times y := \set{z \in \spc{A \times B}: \spc{\pi_A}(z) = x \et \spc{\pi_B}(z) = y},
\]
where $\pi_A : A \times B \to A$ is the projection onto~$A$, and similarly for $\pi_B$.
We shall say that $x$ and $y$ are \intro{\orthogonal} iff $x \times y = \mset{z}$, for some (unique) $z \in \spc{A \times B}$.

We shall say that $x$ and $y$ are \intro{independent} iff for every $z \in x \times y$ we have $\dim z = \dim x + \dim y$.

If $M \elemeq N$, $a \in A(N)$ and $b \in B(N)$, we will say that $b$ and $c$ are \orthogonal over $M$ iff $\tp a M$  and $\tp b M$ are \orthogonal.
If $M \elemeq P$ and $M \elemeq Q$, we shall say that $P$ and $Q$ are \orthogonal over $M$ iff every $b \in P^k$ and $c \in Q^k$ are \orthogonal over $M$, for every $h,k \in \Nat$.
\end{definizione}
What here we called \orthogonal types, are called ``almost orthogonal'' in~\cite{POIZAT:1985}.
We will prove presently that \orthogonal and independent types are the same concept.
\begin{remark}
If $x = \iota a$, for some $a \in A$, then $x$ is \orthogonal to any $M$-type~$y$.
\end{remark}
\begin{lemma}
For every $x \in \spc A$ and $y \in \spc B$, there exists at least one $z \in x \times y$ such that $\dim z = \dim x + \dim y$.
\end{lemma}
\begin{proof}
Assume not.
Then, there would exist $C \subseteq A$ and $D \subseteq B$ definable, such that $x \in \spc C$, $y \in \spc D$, and for every $z \in \spc{C \times D}$, $\dim z < \dim x + \dim y$.
However, this would mean that
\[
\dim C + \dim D = \dim (C \times D) < \dim x + \dim y,
\]
absurd.
\end{proof}
\begin{lemma}
Let $x \in \spc A$ and $y \in \spc B$.
The following are equivalent:
\begin{enumerate}
\item $x$ and $y$ are \orthogonal;
\item $x$ has exactly one extension to $M(y)$;
\item $y$ has exactly one extension to $M(x)$.
\end{enumerate}
\end{lemma}
\begin{lemma}
Let $x$ and $y$ be \orthogonal, with $x' \RKleq x$, and $y' \RKleq y$.
Then, $x'$ and $y'$ are \orthogonal.
\end{lemma}
\begin{proof}
Assume, for contradiction, that $x$ and $y$ are \orthogonal,
but $x'$ and $y'$ are not, where $x'$ and $y'$ are as in the hypothesis.
Therefore,
there exist $z_1 \neq z_2 \in x' \times y'$.
Let $f : M^k \to M^{k'}$ and $g : M^h \to M^{h'}$ be definable, such that $\spf (x) = x'$ and $\spc g(y) = y'$.
Define $h := f \times g :  M^{k+h} \to M^{k' + h'}$.
Let
\[
Z_i := \spc h^{-1}(z_i) \subseteq x \times y, \quad i = 1,2.
\]
Note that $Z_1$ and $Z_2$ are disjoint and non-empty.
However, $x \times y$ is a singleton, which is absurd.
\end{proof}
\begin{lemma}
Let $x \in \spc A$, $y \in \spc B$.
Then, $x$ and $y$ are \orthogonal iff they are independent.
\end{lemma}
\begin{proof}
The ``only if'' direction is immediate from the above Lemma.

For the ``if'' direction, let $z_1 \neq z_2 \in x \times y$.
\Wlog, we can assume that $\dim x = \dim A = h$ and $\dim y = \dim B = k$.
Let $U\subseteq A \times B$ be definable, such that $z_1 \in \spc U$ and $z_2 \notin \spc U$.
Since $x$ and $y$ are independent, $\dim z_1 = \dim z_2 = h+k$.
Hence, \wloG, $U$ is open.
\begin{claim}
$x \times y$ is connected.
\end{claim}
Let $N := M(x)$.
By Lemma~\ref{LEM:TAU}, $x \times y$ is homeomorphic to to $\theta^{-1}(y)$, where $\theta: \Types k {N} \to \Types k M$ is the restriction map.
Hence, by Lemma~\ref{LEM:CONNECTED}, $x \times y$ is connected.

Thus, by Lemma~\ref{LEM:FRONTIER},
$(x \times y) \cap \spc{\fr U} \neq \emptyset$.
Let $z \in (x \times y) \cap \spc{\fr U}$.
Then, $\dim z \leq \dim \fr U < \dim U < h + k$, absurd.
\end{proof}
\begin{example}
If $x \in X \setminus \iota(A)$, then $x$ is \emph{not} \orthogonal to itself.
In fact, let $z(a_1,a_2) \in x \times x$ such that $z$ also satisfies $a_1 = a_2$.
Then, $\dim z = \dim x < \dim x + \dim x$.
\end{example}
\begin{definizione}
We shall say that $M,P,Q,N$ are an \intro{amalgam} iff
\[
M \elemeq P \elemeq N \quad \text{and} \quad
M \elemeq Q \elemeq N.
\]
In this case, we will denote by $PQ$ the elementary sub-structure of $N$ generated by $P \cup Q$.

$M,P,Q,N$ are a \intro{heir-coheir} amalgam iff
for every $\psi(b,c)$ is a formula, and $b \in P^h$, $c \in Q^k$ such that
$N \models \psi(b,c)$, there exists $a \in M^k$ such that
$N \models \psi(b,a)$~\cite{HODGES:1993}.

$M,P,Q,N$ is a \intro{\cotame{}} amalgam iff 
$P/M$ is \tirrational and $Q/M$ is rational.
\end{definizione}
\begin{lemma}
Let $M,P,Q,N$ and $M,P',Q',N'$ be amalgams.
Let $\beta: P \to P'$ and $\gamma: Q \to Q'$ be $M$-isomorphisms.
If $P$ and $Q$ are orthogonal over $M$, then there exists a unique isomorphism $\delta: PQ \to P'Q'$ extending both $\beta$ and $\gamma$.
\end{lemma}
\begin{proof}
Any element of $PQ$ is of the form $f(b,c)$, for some $M$-definable $f: N^{h+k} \to N$, $b \in P^h$, $c \in Q^k$.
Define
\[
\delta\Pa{f(b,c)} := f\Pa{\beta(b), \gamma(c)}.
\]
Since $P$ and $Q$ are orthogonal, $\delta$ is well-defined.
\end{proof}
\begin{lemma}
Let $M,P,Q,N$ be an amalgam.
The following are equivalent:
\begin{enumerate}
\item $M,P,Q,N$ are a heir-coheir amalgam;
\item for every $h,k \in\Nat$, $b \in P^h$, $c \in Q^k$, we have that
$M, M(b), M(c), N$ are a heir-coheir amalgam;
\item for every $h \in \Nat$ and $b \in P^h$,
$\tp b Q$ is a heir of $\tp b M$;
\item for every $k \in \Nat$ and $c \in Q^k$,
$\tp c P$ is a coheir of $\tp c M$;
\item $\tp P Q$ is a heir of $\tp P M$;
\item $\tp Q P$ is a coheir of $\tp Q M$;
\item for every $h,k \in\Nat$, $b \in P^h$, $c \in Q^k$,
$\tp b {M(c)}$ is a heir of $\tp b M$;
\item for every $h,k \in\Nat$, $b \in P^h$, $c \in Q^k$,
$\tp c {M(b)}$ is a coheir of $\tp c M$.
\end{enumerate}
\end{lemma}
\begin{proof}
See \cite{POIZAT:1985}.
\end{proof}
\begin{lemma}\label{LEM:ORTH-HEIR}
Let $M,P,Q,N$ be an amalgam, with $P$ and $Q$ orthogonal over~$M$.
Then, $M,P,Q,N$ is both a heir-coheir and a coheir-heir amalgam.
\end{lemma}
\begin{proof}
By \cite[Theorem~11.01]{POIZAT:1985}, $\tp Q M$ ha at least one heir on~$P$.
Since $Q$ and $P$ are orthogonal, $\tp Q M$ has exactly one extension to~$P$, which therefore must be a heir.
\end{proof}
\begin{lemma}\label{LEM:RAT-HEIR}
Let $M,P,Q,N$ be a coheir-heir amalgam, and $N' := PQ$.
If $Q/M$ is rational, then $N'/P$ is also rational.

Therefore, $\forall d \in N'^n$, such that $d$ is $P$-bounded,
\[
\exists b_0 \in P^n\ d \pequiv b_0, \quad \text{and \emph{a fortiori} } d \mequiv b_0.
\]
\end{lemma}
\begin{proof}
$\tp Q P$ is the (unique) heir of $\tp Q M$.
By \cite[\S11.b]{POIZAT:1985}, $\tp Q P$ is rational, and therefore $N'/P$ is rational.
\end{proof}
\begin{thm}\label{THM:AMALGAM}
Let $M,P,Q,N$ be a \cotame amalgam.
\begin{enumerate}
\item Call $N' := PQ$. 
Then, $N'/Q$ is \tirrational, and $N'/P$ is rational.
\item $M, P, Q, N$ are both a heir-coheir and a coheir-heir amalgam.
\item $P$ and $Q$ are orthogonal over~$M$.
\end{enumerate}
\end{thm}
Proof postponed to \S\ref{SEC:PROOF-AMALGAM}.
\begin{corollary}
Let $x \in \spc A$ be \tirrational, and $y \in \spc B$ be rational.
Then, $x$ and $y$ are \orthogonal.
\end{corollary}
\begin{proof}
By the last point in the Theorem.
\end{proof}
Before proving the Theorem,
we will prove some additional lemmata.
\begin{lemma}\label{LEM:CONTINUITY}
Let $A,B \subseteq M^k$ be \dcompact, and $f:A \to B$ be continuous and definable.
Let $N \extendeq M$, and $c,d \in A(N)$.
If $c \mequiv d$, then $f(c) \mequiv f(d)$.
\end{lemma}
\begin{proof}
Same as Lemma~\cite[1.13]{DRIES-LEW:1995}.
\end{proof}
\begin{example}
In the above lemma, it is essential that $A$ is \dcompact, and $f$ is continuous.
For instance, if $A = \inter{0}{1}^2$ and $\lim_{a \to 0}f(a)$ does not exists (where  $f: A \to [0,1]$ is definable and continuous), then the conclusion does not hold.
\end{example}
\begin{lemma}\label{LEM:1.5}
Let $N/M$ be rational, and $f: N^h \to N^k$ be definable \emph{in $N$}.
Then, the set
\[
D := \set{a \in M^h: f(a) \text{is \mbounded}}
\]
and the map
\begin{align*}
D \to [0,1]^k(M)\\
a   \mapsto \stm\Pa{f(a)}
\end{align*}
are definable in~$M$.
\end{lemma}
\begin{proof}
It is \cite[Corollary~1.5]{DRIES:1997}.
\end{proof}
\begin{lemma}\label{LEM:INTERPOLATE}
Let $M, P,Q,N$ be a \cotame amalgam.
Let $b \in P $ and $c \in C$ such that $b \leq c$.
Then, there exists $a \in M$ such that $b \leq a \leq c$.
If moreover $b \in P \setminus C$, we can also impose $b < a < c$.
\end{lemma}
\subsection{Proof of Theorem~\ref{THM:AMALGAM}}
\label{SEC:PROOF-AMALGAM}
Note that, by Lemma~\ref{LEM:ORTH-HEIR}, the second point is a consequence of the third one.
Moreover, by Lemma~\ref{LEM:RAT-HEIR}, the fact that $N'/P$ is rational is a consequence of the second point.

\begin{paragraph}{First step:} it suffices to consider the case when $P/M$ and $Q/M$ are finitely generated.
\end{paragraph}

The only point that needs clarification is that $N'/Q$ is \tirrational.
Assume therefore that we have proved it for every finite sub-amalgam $M,P',Q',N'$ such that $Q \elemeq Q'$ and $P \elemeq P'$.
We have to prove that every $d\in N'^n$ is $Q$-bounded.

Let $f: M^h \times M^k \to M$ be definable, $b \in P^h$ and $c \in Q^k$ such that $d = f(b,c)$.
%
Since, by assumption, $\tp d {M(c)}$ is \tirrational, $d$ is $M(c)$-bounded, and \emph{a fortiori} $Q$-bounded.

Therefore, we can reduce to the case when $P = M(b)$,  $Q = M(c)$, and $N = M(b,c)$, for some $b \in P^h$, $c \in Q^k$.
\begin{paragraph}{Second step:} it suffice to consider the case when $h =1$, namely $\dim(P/M) = 1$.
\end{paragraph}

We will proceed by induction on~$h$.
Let $b := (b_1, b')$, where $b_1 \in P$ and $b' \in P^{h-1}$.
Let $P_1 := M(b_1)$, and $Q_1 := Q(b_1)$.

Let us apply the case $h = 1$ to the amalgam $M, P_1, Q, Q_1$.
Hence, $Q_1 /P_1$ is rational, and $Q_1 / Q$ is \tirrational.
Moreover, $b_1$ and $c$ are \orthogonal.

Note that $P /P_1$ is \tirrational.
Therefore, $P_1,P,Q_1,N$ are a \cotame amalgam, with $\dim (P /P_1) = h-1$.
Hence, by the inductive case $h-1$, $N'/P$ is rational, and $N'/Q_1$ is \tirrational.
Thus, $N'/Q$ is \tirrational, and the first point is done.


For the last point, let $N''$ be a sufficiently saturated and homogeneous extension of~$N$, and $e \in N''^h$ and $f \in N''^k$ such that $\tp {e} M = \tp b M$ and $\tp {f}{M} = \tp c M$.
We have to prove that $\tp{e,f}{M} = \tp{b,c}M$.
It suffice to show that $\tp{e,c}{M} = \tp{b,c}M$, and $\tp {b,f} M = \tp {b, c}M$.
Namely, that $\tp e Q = \tp b Q$, and $\tp f P = \tp c P$. 

By the case $h = 1$, we have that $\tp {e_1}{Q} = \tp {b_1} Q$, where $e =: (e_1, e')$.
Let $\sigma$ be a $Q$-automorphism of $N''$ such that $\sigma(b_1) = e_1$, and define $g := \sigma(e) = (b_1, g')$.
Thus,
\[
\tp e Q = \tp g Q = \tp {b_1, g'}Q.
\]
Moreover,
\[
\tp {b_1,b'}{M} = \tp b M = \tp e M = \tp g M = \tp {b_1,g'}{M},
\]
and therefore $\tp {b' }{P_1} = \tp {g'}{P_1}$.
Thus, by the case $h - 1$, $\tp {b'} {Q_1} = \tp {g'}{Q_1}$, namely
\[
\tp b {Q} = \tp g {Q} = \tp e {Q}.
\]

Finally, by the case $h = 1$, we have that $\tp f {P_1} = \tp{c}{P_1}$.
Therefore, by the case $h - 1$, we have that $\tp f {P} = \tp c P$.
\begin{paragraph}{Third step:} the case $h = 1$.
\end{paragraph}
\begin{claim}
$(b,c)/M$ are independent.
\end{claim}
If not, then $(b,c)/M$ would be dependent, and therefore $b/M(c)$ is dependent.
Since $c/M$ is independent, this would imply $b \in M(c)$, absurd.


If, for contradiction, $N'/Q$ were not \tirrational, then, since $\dim{N'}/{Q} = 1$, it would be rational.
Hence, $N'/M$ would be rational, and \emph{a fortiori} $P/M$ would be rational, absurd.

For the third point, let $b' \in N$, $c' \in N^k$ such that $\tp {b'}{M} = \tp b M$, and $\tp {c'} M = \tp c M$.
Let $U \subseteq N^{k+1}$ be $M$-definable.
We have to prove that, if $(b,c) \in U$, then $(b',c') \in U$.
\Wlog, we can assume that $U$ is an open cell.
Thus, there exists an open cell $W \subseteq N^k$ and functions
$f,g : W \to N$, all definable (over~$M$) such that $f < g $ and
\[
U = \set{(t,u)\in N^{k+1}: u \in W, f(u)< t < g(u)}.
\]
Thus, $f(c) < b < g(c)$.
If $f(c) = a_1 \in M$, then $f(c) = f(c') = a_1 < b'$, and similarly for~$g$.
Otherwise, there exist $a_1, a_2 \in M$ such that
\[
f(c) < a_1 < b < a_2 < g(c).
\]
Thus,
\[
f(c') < a_1 < b' < a_2 < g(c'),
\]
and $(b',c') \in U$.
\qed
\begin{example}
Let $M$ be the field of real algebraic numbers, $b_1, b_2 \in \Real$ be algebraically independent, and $\tp{t}{\Real} = 0^+$.
Let $P := M(b_1,b_2)$, $Q := M(t, b_1 + t b_2)$, $N := PQ = M(b_1,b_2,t)$.
Then, $P \cap Q = M$, but $P$ and $Q$ are \emph{not} orthogonal over $M$, because $\dim N/M = 3 < \dim P/m + \dim Q/M$.
\end{example}




%% file: main.bbl
\def\cprime{$'$}
\begin{thebibliography}{10}

\bibitem{BAI-POIZAT:1998}
Y.~Baisalov and B.~Poizat.
\newblock Paires de structures o-minimales.
\newblock {\em J. Symbolic Logic}, 63(2):570--578, 1998.

\bibitem{COSTE-CARRAL:1983}
M.~Carral and M.~Coste.
\newblock Normal spectral spaces and their dimensions.
\newblock {\em J. Pure Appl. Algebra}, 30(3):227--235, 1983.

\bibitem{DELFS:1985}
H.~Delfs.
\newblock The homotopy axiom in semialgebraic cohomology.
\newblock {\em J. Reine Angew. Math.}, 355:108--128, 1985.

\bibitem{HODGES:1993}
W.~Hodges.
\newblock {\em Model theory}, volume~42 of {\em Encyclopedia of Mathematics and
  its Applications}.
\newblock Cambridge University Press, Cambridge, 1993.

\bibitem{JONES:PHD}
G.~O. Jones.
\newblock {\em Local to global methods in o-minimal expansions of fields}.
\newblock PhD thesis, University of Oxford, Wolfson College, 2006.

\bibitem{PETERZIL-STARCHENKO:1998}
Y.~Peterzil and S.~Starchenko.
\newblock A trichotomy theorem for o-minimal structures.
\newblock {\em Proc. London Math. Soc. (3)}, 77(3):481--523, 1998.

\bibitem{PILLAY:1987}
A.~Pillay.
\newblock First order topological structures and theories.
\newblock {\em J. Symbolic Logic}, 52(3):763--778, 1987.

\bibitem{PILLAY:1988}
A.~Pillay.
\newblock Sheaves of continuous definable functions.
\newblock {\em J. Symbolic Logic}, 53(4):1165--1169, 1988.

\bibitem{PILLAY:1994}
A.~Pillay.
\newblock Definability of types, and pairs of {O}-minimal structures.
\newblock {\em J. Symbolic Logic}, 59(4):1400--1409, 1994.

\bibitem{POIZAT:1985}
B.~Poizat.
\newblock {\em Cours de th{\'e}orie des mod{\`e}les}.
\newblock Bruno Poizat, Lyon, 1985.
\newblock Une introduction {\`a} la logique math{\'e}matique contemporaine. [An
  introduction to contemporary mathematical logic].

\bibitem{DRIES:1997}
L.~van~den Dries.
\newblock {$T$}-convexity and tame extensions. {II}.
\newblock {\em J. Symbolic Logic}, 62(1):14--34, 1997.

\bibitem{DRIES:1998}
L.~van~den Dries.
\newblock {\em Tame topology and o-minimal structures}, volume 248 of {\em
  London Math. Soc. Lecture Note Ser.}
\newblock Cambridge University Press, Cambridge, 1998.

\bibitem{DRIES-LEW:1995}
L.~van~den Dries and A.~H. Lewenberg.
\newblock {$T$}-convexity and tame extensions.
\newblock {\em J. Symbolic Logic}, 60(1):74--102, 1995.

\bibitem{WILKIE:2003}
A.~J. Wilkie.
\newblock Covering definable open sets by open cells.
\newblock In M.~Edmundo, D.~Richardson, and A.~Wilkie, editors, {\em O-minimal
  Structures, Proceedings of the RAAG Summer School Lisbon 2003}, Lecture Notes
  in Real Algebraic and Analytic Geometry. Cuvillier Verlag, 2005.

\end{thebibliography}
